\documentclass[11pt]{article}
\usepackage{amsmath}
\usepackage{graphicx}
\usepackage{amsfonts}
\usepackage[catalan,spanish,english]{babel} 
\usepackage[psamsfonts]{amssymb} 
\usepackage{cmmib57} 
\usepackage{exscale} 
\usepackage{amscd} 
\usepackage{latexsym} 
\usepackage{graphicx} 
\allowdisplaybreaks
\numberwithin{equation}{section}
\newtheorem{stat}{Statement}[section]
\newtheorem{theorem}[stat]{Theorem}

\newtheorem{proposition}[stat]{Proposition}
\newtheorem{lemma}[stat]{Lemma}
\newtheorem{definition}[stat]{Definition}

\newcommand{\beq}{\begin{equation}}
\newcommand{\eeq}{\end{equation}}
\newcommand{\beqn}{\begin{equation*}}
\newcommand{\eeqn}{\end{equation*}}
\begin{document}
\begin{titlepage}
\medskip

\begin{center}
{\Large\bf Mild Solutions for a Class 
of Fractional SPDEs
and Their Sample Paths}
\smallskip

by\\
\vspace{7mm}
\begin{tabular}{l@{\hspace{10mm}}l@{\hspace{10mm}}l}
{\sc Marta Sanz-Sol\'e}$\,^{(\ast)}$ &and& {\sc Pierre-A. Vuillermot}\\
{\small Facultat de Matem\`atiques} && {\small UMR-CNRS 7502}\\
{\small Universitat de Barcelona} && {\small Institut \'Elie Cartan} \\
{\small Gran Via 585} && {\small BP 239}\\
{\small E-08007 Barcelona, Spain} && {\small F-54506 Vandoeuvre-l\`es Nancy Cedex }\\
{\small marta.sanz@ub.edu }&&{\small vuillerm@iecn.u-nancy.fr}\\
\null
\end{tabular}
\end{center}
\vskip 1.5cm

\noindent{\bf Abstract.} We introduce a notion of mild solution for a class of
non-autonomous parabolic stochastic partial differential equations 
defined on a bounded open subset $D\subset\mathbb{R}^{d}$ and driven by an
infinite-dimensional fractional noise. 
We prove the existence of such a solution, establish its relation with the variational solution introduced in \cite{nuavu}
and  the H\"{o}lder continuity of its sample paths when we consider it as an
$L^{2}(D)$--valued stochastic process. When $h$ is an affine function,
we also prove uniqueness.
 An immediate consequence
of our results is the  indistinguishability of mild and variational
solutions in the case of uniqueness.
 \bigskip

\smallskip

\noindent{\bf AMS Subject Classification}\\
{\sl Primary}: 60H15, 35R60, 35K55.\\
{\sl Secondary}: 60G17, 60G15, 60G18.
\medskip

\noindent
\footnotesize
{\begin{itemize} \item[$^{(\ast)}$] Supported by the grant MTM 2006-01351 from the \textit{Direcci\'on General de
Investigaci\'on, Ministerio de Educaci\'on y Ciencia, Spain.}
\end{itemize}}
\end{titlepage}

\newpage

\section{Introduction and Outline}
\label{s1}
In the last decades, the interest in fractional Brownian motion, first introduced in \cite{kolmogorov} and referred to as {\sl fBm} in the sequel, has increased enormously 
as one important ingredient of fractal models in the sciences. The paper \cite{mandel} has been one of the
keystones that has attracted the attention of part of  the probabilistic community to this challenging object.
Some of the research on {\sl fBm} has significantly influenced the present state of the art of Gaussian processes
(see for instance  \cite{bn}, \cite{np}, \cite{pipitaqu}, just to mention 
a few). An important aspect of the study of {\sl fBm} lies in the domain of stochastic analysis. Since this process is neither a semimartingale 
nor a Markov process, It\^o's theory does not apply. For values of the Hurst parameter $H$ greater than $\frac{1}{2}$ -the regular case- integrals of Young's
type and fractional calculus techniques have been considered (\cite{zaehle1}, \cite{zaehle2}).
However, for $H$ less than $\frac{1}{2}$ this approach fails. The integral representation of {\sl fBm}
as a Volterra integral with respect to the standard Brownian motion has been successfully exploited in setting up a stochastic calculus where classical tools of 
Gaussian processes along with fractional and Malliavin calculus are combined.
A pioneering work in this context is \cite{decreustu}, then \cite{alomanu2}, \cite{cacomon},
\cite{d} and also \cite{feypradel}.
Since then, there have been many contributions to the subject. Let us refer to \cite{n} for enlightening contents
and a pretty complete list of references. Rough path analysis (see  \cite{lq}) provides a new approach somehow related to Young's
approach. 

One of the main reasons for developing a stochastic calculus based on {\sl fBm} is mathematical modeling. The theory of ordinary and partial differential equations
driven by a fractional noise is nowadays a very active field of research. Some of the motivations come from a number of applications in engineering, biophysics and 
mathematical finance; to refer only to a few, let us mention \cite{dms}, \cite{ks}, \cite {roger}. There are also purely mathematical motivations. Problems studied so far 
range from the existence, the uniqueness, the regularity and the long-time behaviour of solutions to large deviations, support theorems and the analysis of the law of the solutions using Malliavin calculus. 
Without aiming to be exhaustive, let us refer to \cite{b-h}, \cite{f},  \cite{guletin}, \cite{huokszang}, 
 \cite{klinzaehle}, \cite{lq},
\cite{maslonua}, \cite{m-ss}, \cite{nuarasca}, \cite{nuavu}, \cite{q-t} and \cite{tintuvien}) 
for a reduced sample of published work. 

This paper aims to pursue the investigations of \cite{nuavu}, where 
the authors develop an existence and uniqueness theory of variational
solutions for a class of non-autonomous semilinear partial differential
equations driven by an infinite-dimensional multiplicative fractional noise
through the construction and the convergence of a suitable Faedo-Galerkin scheme.

As is the case for deterministic partial differential equations, a recurrent
difficulty  is the necessity to decide \textit{ab
initio} what solution concept is relevant, since there are several \textit{a
priori} non-equivalent possibilities to choose from. Thus, while in
 \cite{nuavu}  two notions of \textit{variational solution} that
are subsequently proved to be indistinguishable are introduced, the focus in \cite{guletin}
or \cite{maslonua} is rather on the idea of\textit{\ mild solution}, 
namely, a solution which can be expressed as a nonlinear integral
equation that involves the linear propagator of the theory without any
reference to specific classes of test functions. Consequently, this leaves
entirely open the question of knowing whether the variational and mild notions
are in some sense equivalent, and indeed we are not aware of any connections
between them thus far in this context. For equations of the type considered in this
article but driven by standard Wiener processes, this issue was addressed
in \cite{sansovu1}. In \cite{denistoica} a similar question was analyzed for a class of very  general 
SPDEs driven by a finite-dimensional Brownian motion.

In this article we consider the same class of equations as in \cite{nuavu}. 
We develop an existence and uniqueness  theory of mild solutions and  prove the indistinguishability of variational and mild solutions.
We also prove the H\"older continuity of their sample paths.

Before defining the class of problems we shall investigate, let us fix the notation.
All the functional spaces we introduce are real and we use the standard
notations for the usual Banach spaces of differentiable functions, of
H\"{o}lder continuous functions, of Lebesgue integrable functions and for the
related scales of Sobolev spaces defined on regions of Euclidean space used for instance in  \cite{adams}. 
For $d\in\mathbb{N}^+$ let $D\subset\mathbb{R}^{d}$ be an open and  bounded set whose boundary $\partial D$ is
of class $\mathcal{C}^{2+\beta}$ for some $\beta\in(0,1)$ (see, for instance,
\cite{eidelzhita} and \cite{ladyuralsolo} for a definition of
this and related concepts). We will denote by $(.,.)_2$ the standard inner product in 
$L^2({\mathbb{R}}^d)$, by $(.,.)_{\mathbb{R}^d}$ the Euclidean inner product in $\mathbb{R}^d$
and by $|.|$ the associated Euclidean norm.

 Let $(\lambda_{i})_{i\in\mathbb{N}^{+}}$ be any sequence of positive real numbers such that
$\sum_{i=1}^{+\infty}\lambda_{i}<+\infty$ and $(e_{i})_{i\in\mathbb{N}^{+}}$  an orthonormal basis of $L^{2}(D)$ such that  
$\sup_{i\in\mathbb{N}^{+}}\left\|  e_{i}\right\|  _{\infty}<+\infty$ (the existence of such a basis follows from  
\cite{ovespelc}). We then define the linear, self-adjoint, positive,
non-degenerate trace-class operator $C$ in $L^{2}(D)$ by $Ce_{i}=\lambda
_{i}e_{i}$ for each $i$. In the sequel we write $\left(  \left(  B_{i}
^{H}(t)\right)  _{t\in\mathbb{R}^{+}}\right)  _{i\in
\mathbb{N}^{+}}$ for a sequence of one-dimensional, independent, identically
distributed fractional Brownian motions with Hurst parameter 
$H\in(0,1)$, defined on the complete probability space
$\left(  \Omega,\mathcal{F},\mathbb{P}\right)  $ and starting at the origin.
We introduce the $L^{2}(D)$-valued fractional Wiener process $\left(
W^{H}(.,t)\right)  _{t\in\mathbb{R}^{+}}$ by setting
\begin{equation}
W^{H}(.,t):=\sum_{i=1}^{+\infty}\lambda_{i}^{\frac{1}{2}}
e_{i}(.)B_{i}^{H}(t), \label{1}
\end{equation}
where the series converges a.s. in the strong topology of $L^{2}(D)$, by virtue
of the basic properties of the $B_{i}^{H}(t)$'s and the fact
that $C$ is trace-class. 

Let $T>0$, $\alpha\in(1-H,\frac{1}{2})$ and $(F(t), t\in [0,T] )$ be a stochastic process taking values in the space of 
linear bounded operators on $L^2(D)$. Assume that
\beqn
\sup_{i\in\mathbb{N}^+}\Vert F(s)e_i\Vert_{\alpha,1}<+\infty,
\eeqn
where for a function $f:[0,T]\to L^2(D)$,
\beqn
\Vert f\Vert_{\alpha,1}= \int_0^T\left(\frac{\Vert f(s)\Vert_2}{s^\alpha}+\int_0^s \frac{\Vert f(s)-f(r)\Vert_2}{(s-r)^{\alpha+1}}dr \right) ds.
\eeqn
Following \cite{maslonua} we can define a pathwise generalized Stieltjes integral (see also \cite{zaehle1})
\beqn
\int_0^T F(s) W^H(ds):=\sum_{i=1}^{+\infty} \lambda_i^{\frac{1}{2}} \int_0^T F(s)e_i B_i^H(ds),
\eeqn
which satisfies the property
\begin{align*}
\left\Vert\int_0^T F(s) W^H(ds)\right\Vert_2&\le \sup_{i\in\mathbb{N}^+}\Vert F(s)e_i\Vert_{\alpha,1}
\left(\sum_{i=1}^{+\infty} \lambda_i^{\frac{1}{2}}\Lambda_\alpha(B_i^H)\right).
\end{align*}
Here $\Lambda_\alpha(B_i^H)$ is a positive random variable defined in terms of a Weyl derivative (see \cite{maslonua}, Equation (2.4)),
satisfying $\sup_{i\in\mathbb{N}^+}E\left(\Lambda_\alpha(B_i^H)\right)<+\infty$, as is proved in \cite{nuarasca}, Lemma 7.5.
Consequently, if $\sum_{i=1}^{+\infty} \lambda_i^{\frac{1}{2}}<+\infty$, the random variable 
\begin{equation}
r_{\alpha}^{H}:=\sum_{i=1}^{+\infty}\lambda_{i}^{\frac{1}{2}}\Lambda_{\alpha}(B_{i}^{H}). \label{29}
\end{equation}
is finite, a.s., and then
\beq
\label{i}
\left\Vert\sum_{i=1}^{+\infty} \lambda_i^{\frac{1}{2}} \int_0^T F(s)e_i B_i^H(ds)\right\Vert_2\le 
r_{\alpha}^H \sup_{i\in\mathbb{N}^+}\Vert F(s)e_i\Vert_{\alpha,1}.
\eeq
\smallskip

Next, we introduce the class of real, parabolic,
initial-boundary value problems formally given by
\begin{align}
du(x,t)  &  = \left(\operatorname*{div}(k(x,t)\nabla
u(x,t))+g(u(x,t))\right)  dt+h(u(x,t))W^{H}(x,dt),\nonumber\\
(x,t)  &  \in D\times \left(  0,T\right]  ,\nonumber\\
u(x,0)  &  =\varphi(x),\text{ \ \ }x\in\overline{D},\nonumber\\
\frac{\partial u(x,t)}{\partial n(k)}  &  =0,\text{ \ \ }(x,t)\in\partial
D\times\left(  0,T\right]  , \label{2}
\end{align}
where the last relation stands for the conormal derivative of $u$ relative to the matrix-valued field $k$.

In the next section we shall give a rigorous meaning to such a  formal expression and for this, we shall use the 
pathwise integral described before.

In the sequel we write $n(x)$ for the unit outer normal vector at $x\in$ $\partial D$ and introduce 
the following set of assumptions:
\begin{description}
\item{(C)} The square root $C^{\frac{1}{2}}$ of the covariance operator is
trace-class, that is, we have $\sum_{i=1}^{+\infty}\lambda_{i}^{\frac{1}{2}}<+\infty$.
\item{(${\rm K}_{\beta,\beta^{\prime}}$)} The entries of $k$ satisfy $k_{i,j}%
(.)=k_{j,i}(.)$ for all $i,j\in\left\{  1,...,d\right\}  $ and there exists a
constant $\beta^{\prime}\in\left(  \frac{1}{2},1\right]  $ such that
$k_{i,j}\in\mathcal{C}^{\beta,\beta^{\prime}}(\overline{D}\times\left[0,T\right])$
for each $i,j$. In addition, we have
$k_{i,j,x_{l}}:=\frac{\partial k_{i,j}}{\partial x_{l}}\in\mathcal{C}
^{\beta,\frac{\beta}{2}}(\overline{D}\times\left[  0,T\right])$
for each $i,j,l$ and there exists a constant $\underline{k}\in\mathbb{R}%
_{\ast}^{+}$ such that the inequality
$
\left(  k(x,t)q,q\right)  _{\mathbb{R}^{d}}\geq\underline{k}\left|  q\right|
^{2}
$
holds for all $q\in\mathbb{R}^{d}$ and all $(x,t)\in\overline{D}\times\left[
0,T\right]$.\\
Finally, we
have
$$
(x,t)\mapsto\sum_{i=1}^{d}k_{i,j}(x,t)n_{i}(x)\in\mathcal{C}^{1+\beta
,\frac{1+\beta}{2}}(\partial D\times\left[  0,T\right])
$$
for each $j$ and the conormal vector-field $(x,t)\mapsto
n(k)(x,t):=k(x,t)n(x)$ is outward pointing, nowhere tangent to $\partial D$
for every $t$.
\item{(L)} The functions $g,h:\mathbb{R\mapsto R}$ are Lipschitz continuous.
\item{(I)} The initial condition satisfies $\varphi\in\mathcal{C}^{2+\beta}%
(\overline{D})$ and the conormal boundary condition relative to $k$.
\end{description}



Finally, we consider the following assumption which also appears in \cite{nuavu}: 
\smallskip

\noindent (${\rm H}_\gamma$)   The derivative $h^{\prime}$ is H\"{o}lder continuous with exponent $\gamma\in\left(0,1\right]$
and bounded; moreover, the Hurst parameter satisfies $H\in\left(\frac{1}{\gamma+1},1\right)$.

Notice that if the derivative $h^{\prime}$ is Lipschitz continuous this amounts to assuming $H\in\left(\frac{1}{2},1\right)$.

Problem (\ref{2}) is identical to the
initial-boundary value problem investigated in \cite{nuavu}, up to Hypotheses
(K$_{\beta,\beta^{\prime}}$) which imply Hypotheses (K)
of that article. This
immediately entails the existence of what is called there a variational
solution of type II for (\ref{2}), henceforth simply coined variational solution. 
With (K$_{\beta,\beta^{\prime}}$) we have the existence and regularity properties of the Green function
associated with the differential operator governing (\ref{2}). We shall give more details on this in the 
next section.
\medskip

We organize this article in the following way. In Section \ref{s2} we first recall
the notion of {\it variational solution} and introduce a
notion of {\it mild solution} for (\ref{2}) by means of a family of evolution
operators in $L^{2}(D)$ generated by the corresponding deterministic 
Green's function. We then proceed by stating our main
results concerning the existence, uniqueness, and H\"{o}lder regularity of the mild
solution along with its indistinguishability from the
variational solution when $h$ is an affine function. 
The section ends with a discussion about the results and methods of their proofs. 
These are gathered in Section \ref{s3}.
\bigskip


\section{Statement and Discussion of the Results}
\label{s2}

In the remaining part of this article we write $H^{1}(D\times(0,T))$ for the
isotropic Sobolev space on the cylinder $D\times(0,T)$, which consists of all
functions $v\in L^{2}(D\times(0,T))$ that possess distributional derivatives
$v_{x_{i}}$, $v_{\tau}\in L^{2}(D\times(0,T))$. 
The set of all $v\in H^{1}(D\times(0,T))$ which do not depend
on the time variable identifies with $H^{1}(D)$, the usual Sobolev
space on $D$ whose 
norm we denote by 
$\left\|  .\right\|  _{1,2}$. 

For $0<\alpha<1$
we introduce the Banach space $\mathcal{B}^{\alpha,2}(0,T;L^{2}(D))$
of all Lebesgue-measurable mappings $u:\left[  0,T\right]  \mapsto L^{2}(D)$
endowed with the norm 
\begin{equation}
\left\|  u\right\|  _{\alpha,2,T}^{2}:= \left(\sup_{t\in\left[0,T\right]}\left\|u(t)\right\|_{2}\right)^2 + \int_{0}^{T}dt\left(  \int_{0}
^{t}d\tau\frac{\left\|  u(t)-u(\tau)\right\|  _{2}^{{}}}{(t-\tau)^{\alpha+1}
}\right)  ^{2}<+\infty. \label{3}
\end{equation}
Notice that $\Vert \cdot\Vert_{\alpha,1}\le c\left\| \cdot \right\|  _{\alpha,2,T}$,  
and also that the spaces $\mathcal{B}^{\alpha,2}(0,T;L^{2}(D))$ decrease when $\alpha$ increases. 

We recall the following notion introduced in \cite{nuavu}, in which the function $x\mapsto v(x,t)\in L^{2}(D)$
is interpreted as the Sobolev trace of $v\in H^{1}(D\times(0,T))$ on the
corresponding hyperplane. 

\begin{definition}
\label{d1}
Fix $H\in\left(\frac{1}{2},1\right)$ and let $\alpha\in\left(1-H,\frac{1}{2}\right)$. We
assume that conditions (C), (L) are satisfied and that the initial condition $\varphi$ belongs to $L^2(D)$. In addition we suppose that  the symmetric matrix valued function $k$ satisfies
\beqn
\underline k |q|^2 \le (k(x,t)q,q)_{\mathbb{R}^d}\le \overline k |q|^2,
\eeqn
for any $q\in\mathbb{R}^d$ and some positive constants $\underline k$, $\overline k$
independent of $x$ and $t$.

Under these conditions, the $L^{2}(D)$-valued random field $\left( 
u_{V}(.,t)\right)  _{t\in\left[  0,T\right]  }$ defined and measurable on
$\left(  \Omega,\mathcal{F},\mathbb{P}\right)  $ is a \textit{variational
solution }to Problem (\ref{2}) if:

(1)  $u_{V}\in L^{2}(0,T;H^{1}(D))\cap\mathcal{B}^{\alpha,2}%
(0,T;L^{2}(D))$ a.s., which means that 
\[
\int_{0}^{T}dt\left\|  u_{V}(.,t)\right\|  _{1,2}^{2}=\int_{0}^{T}dt\left(
\left\|  u_{V}(.,t)\right\|  _{2}^{2}+\left\|  \nabla u_{V}(.,t)\right\|
_{2}^{2}\right)  <+\infty
\]
and $\left\|  u_{V}\right\|  _{\alpha,2,T}<+\infty$
hold a.s.

(2) The integral relation
\begin{align}
\int_{D}dx\ v(x,t)u_{V}(x,t)  &  =\int_{D}dx\ v(x,0)\varphi(x)+\int_{0}^{t}
d\tau\int_{D}dx\ v_{\tau}(x,\tau)u_{V}(x,\tau)\nonumber\\
&  -\int_{0}^{t}d\tau\int_{D}dx\left(  \nabla v(x,\tau),k(x,\tau)\nabla
u_{V}(x,\tau)\right)  _{\mathbb{R}^{d}}\nonumber\\
&  +\int_{0}^{t}d\tau\int_{D}dxv(x,\tau)g(u_{V}(x,\tau))\nonumber\\
&+\sum_{i=1}^{+\infty}\lambda_{i}^{\frac{1}{2}}\int_{0}^{t}\left(v(.,\tau),h(u_{V}(.,\tau))e_{i}\right)  _{2}B_{i}^{H}(d\tau).
 \label{4}
\end{align}
holds a.s. for every $v\in$ $H^{1}(D\times(0,T))$ and every $t\in\left[
0,T\right] $.
\end{definition}


With the standing hypotheses we easily infer that each term in
(\ref{4}) is finite a.s. In particular, Hypothesis (C) and the fact that $h$ is Lipschitz continuous,
along with (\ref{i})  imply the
absolute convergence, a.s., of the series of the last term in (\ref{4}).


\bigskip

Let
$G:\overline{D}\times\left[  0,T\right]  \times\overline{D}\times\left[  0,T\right]
\diagdown\left\{  s,t\in\left[  0,T\right]  :s\geq t\right\} \to\mathbb{R}$
be the parabolic Green's function associated with
the principal part of (\ref{2}). Assume that (K$_{\beta,\beta^{\prime}}$) holds; it is well-known that 
$G$ is a continuous function, twice continuously differentiable in $x$, once continuously
differentiable in $t$. For every $(y,s)\in D\times (0,T]$, it is also a classical solution to the linear initial-boundary
value problem
\begin{align}
\partial_{t}G(x,t;y,s)  &  =\operatorname*{div}(k(x,t)\nabla_{x}
G(x,t;y,s)),\text{ \ \ }(x,t)\in D\times\left(0,T\right],\nonumber\\
\frac{\partial G(x,t;y,s)}{\partial n(k)}  &  =0,\text{ \ \ \ }(x,t)\in
\partial D\times\left(  0,T\right], \label{6}
\end{align}
with
\[
\int_{D}dyG(.,s;y,s)\varphi(y):=\lim_{t\searrow s}\int_{D}dyG(.,t;y,s)\varphi
(y)=\varphi(.),
\]
and satisfies the heat kernel estimates
\begin{equation}
\left|  \partial_{x}^{\mu}\partial_{t}^{\nu}G(x,t;y,s)\right|  \leq
c(t-s)^{-\frac{d+\left|  \mu\right|  +2\nu}{2}}\exp\left[  -c\frac{\left|
x-y\right|  _{{}}^{2}}{t-s}\right]  \label{7}
\end{equation}
for $\mu=(\mu_{1},...,\mu_{d})\in\mathbb{N}^{d}$, $\nu\in\mathbb{N}$ and
$\left|\mu\right| +2\nu\leq 2$, with $\left|  \mu\right|  =\sum_{j=1}^{d}
\mu_{j}$ (see, for instance, 
\cite{eidelzhita} or
\cite{ladyuralsolo}). In particular, for $\left|  \mu\right|  =\nu=0$ we have
\begin{equation}
\left|  G(x,t;y,s)\right|  \leq c(t-s)^{-\frac{d}{2}}\exp\left[
-c\frac{\left|  x-y\right|  _{{}}^{2}}{t-s}\right]. \label{8}
\end{equation}
We shall refer to (\ref{8}) as the {\it Gaussian
property} of $G$. 

We can now define the notion of
mild solution for (\ref{2}).

\begin{definition}
\label{d2}
Fix $H\in\left(\frac{1}{2},1\right)$ and let $\alpha\in\left(1-H,\frac{1}{2}\right)$. Assume that the hypotheses (C), (K$_{\beta,\beta^{\prime}}$), (L)
hold and that the initial condition $\varphi$ is bounded.

Under these assumptions, the $L^{2}(D)$-valued random field $\left(
u_{M}(.,t)\right)  _{t\in\left[  0,T\right]  }$ defined and measurable on
$\left(  \Omega,\mathcal{F},\mathbb{P}\right)  $ is a \textit{mild}
\textit{solution }to Problem (\ref{2}) if the following two conditions are satisfied:

(1) $u_{M}\in L^{2}(0,T;H^{1}(D))\cap\mathcal{B}^{\alpha,2}
(0,T;L^{2}(D))$ a.s.

(2) The relation
\begin{align}
u_{M}(.,t)  &  =\int_{D}dy\ G(.,t;y,0)\varphi(y)+\int_{0}^{t}d\tau\int
_{D}dy\ G(.,t;y,\tau)g\left(  u_{M}(y,\tau)\right) \nonumber\\
&+\sum_{i=1}^{+\infty}\lambda_{i}^{\frac{1}{2}}\int_{0}^{t}\left(  \int
_{D}dy\ G(.,t;y,\tau)h\left(u_{M}(y,\tau)\right)  e_{i}(y)\right)
B_{i}^{H}(d\tau) \label{9}
\end{align}
holds a.s. for every $t\in\left[  0,T\right]  $ as an equality in $L^{2}(D)$.
\end{definition}

We shall prove in Lemma \ref{l2} that with the standing assumptions, each term in (\ref{9}) indeed
defines a $L^{2}(D)$-valued stochastic process.
\medskip

The main results of this article are gathered in the next theorem.
\bigskip

\begin{theorem}
\label{t1}
Assume that Hypotheses ($C$), ($K_{\beta,\beta^{\prime}}$), ($L$), ($I$) and ($H_{\gamma}$) hold; then
the following statements are valid:

\begin{description}
\item{(a)} Fix $H\in\left(\frac{1}{\gamma+1},1\right)$ and let $\alpha\in\left(1-H,\frac{\gamma}{\gamma+1}\right)$. Then
Problem (\ref{2}) possesses a variational solution
$u_V$; moreover, every such variational solution is a mild solution $u_M$ to (\ref{2}).
More presicely, for every $t\in\left[  0,T\right]$, $u_V(.,t)=u_M(.,t)$  a.s. in $L^{2}(D)$.


\item{(b)} Fix $H\in\left(\frac{1}{\gamma+1}\vee\frac{d+1}{d+2},1\right)$ and then $\alpha\in\left(1-H,\frac{\gamma}{\gamma+1}\wedge\frac{1}{d+2}\right)$. Assume in addition that $h$ is an affine function. Then $u_V$  is the unique variational solution to
 (\ref{2}), while $u_M$ is its unique mild solution. 
 
\item{(c)} Let $H$ and $\alpha$ be as in part (b).
Then every mild solution $u_{M}$ to Problem (\ref{2}) 
is H\"{o}lder continuous with respect to the time variable. More precisely,
there exists a positive, a.s. finite  random variable $R_{\alpha}^{H}$
such that the estimate
\begin{equation}
\left\|  u_{M}(.,t)-u_{M}(.,s)\right\|  _{2}\leq R_{\alpha}^{H}
\left|  t-s\right|  ^{\theta}\left(  1+\left\|u_{M}\right\|  _{\alpha,2,T}\right)
\label{h1}
\end{equation}
holds a.s. for all $s,t\in\left[  0,T\right] $ and every
$\theta\in\left( 0,\left(\frac{1}{2}-\alpha\right)\wedge\frac{\beta}{2}\right)$.
 \end{description}
\end{theorem}



{\bf Remarks}
\begin{enumerate}
\item The existence of a mild solution will be proved by reference to
the existence of a variational solution. 
This is in contrast with the method of \cite{maslonua}, in which the
authors prove the existence of mild solutions for a class of
\textit{autonomous}, parabolic, fractional stochastic initial-boundary value
problems by means of Schauder's fixed point theorem. Their method thus
requires the construction of a continuous map operating in a compact and
convex set of a suitable functional space.   
If $h$ is an affine function, the arguments of the proof of Statement (b)
(see (\ref{322})) show that a similar approach might be possible for
our equation.
To the best of our knowledge, there exists as yet no such direct way to prove the
existence of mild solutions to (\ref{2}) for a non affine $h$. 

\item If $h$ is not affine, the question of uniqueness remains unsettled. 
In fact, uniqueness could be proved if we were able to extend the  inequality (\ref{322})
to any Lipschitz function $h$. This does not seem to be a trivial point, due to the form of the second term in the right-hand side of (\ref{3}).

\item If $h$ is an affine function, Theorem \ref{t1} establishes the complete
indistinguishability of mild and variational solutions, although we do not
know whether this property still holds for a general $h$.

\item If $h$ is a constant function the setting of the problems and their proofs become much 
simpler. Indeed, in Definitions \ref{d1} and \ref{d2} the space $\mathcal{B}^{\alpha,2}(0,T; L^2(D))$
can be replaced by the larger one $L^\infty(0,T; L^2(D))$, consisting of Lebesgue-measurable
mappings $u:[0,T] \to L^2(D)$ such that $\sup_{t\in[0,T]}\Vert u(t)\Vert_2 <\infty$. This can be checked
by going through the proofs of  \cite{nuavu} and Lemma \ref{l2} of Section \ref{s3}.
Moreover, the range of values of $\theta$ in statement (c) of Theorem \ref{t1} can be extended to the interval $\left(0,\frac{\beta}{2}\right)$.
This can be easily checked by going through the proof of Proposition \ref{p4}, by checking first that the right-hand side
of the inequalities (\ref{40}), (\ref{41}) can be replaced by $c(t-s)^{\delta} (s-\tau)^{-\delta}$ and 
$c(t-s)^{\frac{\delta}{2}} (s-\tau)^{-\frac{1}{2}}(\tau-\sigma)^{\frac{1}{2}(1-\delta)}$, respectively.

\item 
By using the \textit{factorization method}, and under a different set of assumptions on the range of admissible values of $H$ and $\alpha$,
we can obtain a different range of values for the H\"older exponent which in general do not provide as good an estimate as (\ref{h1}) does.
We deal with this question in Proposition \ref{pfact}. 
The \textit{factorization method}
has been introduced in \cite{DP-K-Z} and since then
extensively used for the analysis of the sample paths of solutions to
parabolic stochastic partial differential equations (see, for instance,
\cite{sansovu1}). 
\end{enumerate}

\section{Proofs of the Results}
\label{s3}

In what follows we write $c$ for all the irrelevant deterministic constants that occur in
the various estimates. 
We begin by
recalling that the uniformly elliptic partial differential operator with
conormal boundary conditions in the principal part of (\ref{2}) admits a
self-adjoint, positive realization $A(t):=-\operatorname*{div}(k(.,t)\nabla)$
in $L^{2}(D)$ on the domain
\begin{equation}
\mathcal{D}(A(t))=\left\{  v\in H^{2}(D):\left(  \nabla
v(x),k(x,t)n(x)\right)  _{\mathbb{R}^{d}}=0,\text{ \ }(x,t)\in\partial
D\times\left[  0,T\right]  \right\}  \label{11}
\end{equation}
(see, for instance, \cite{lions}). An important consequence
of this property is that the parabolic Green's function $G$ is also, for every
$(x,t)\in D\times\left(  0,T\right]  $ with $t>s$, a classical solution to the
linear boundary value problem
\begin{align}
\partial_{s}G(x,t;y,s)  &  =-\operatorname*{div}(k(y,s)\nabla_{y}
G(x,t;y,s)),\text{ \ \ }(y,s)\in D\times\left(  0,T\right]  ,\nonumber\\
\frac{\partial G(x,t;y,s)}{\partial n(k)}  &  =0,\text{ \ \ \ }(y,s)\in
\partial D\times\left(  0,T\right]  , \label{12}
\end{align}
dual to (\ref{6}) (see, for instance, \cite{eidelzhita} or \cite{friedman});
this means that along with (\ref{7}) we also have 
\begin{equation}
\left|  \partial_{y}^{\mu}\partial_{s}^{\nu}G(x,t;y,s)\right|  \leq
c(t-s)^{-\frac{d+\left|  \mu\right|  +2\nu}{2}}\exp\left[  -c\frac{\left|
x-y\right|  _{{}}^{2}}{t-s}\right]  \label{13}
\end{equation} 
for $\left|  \mu\right|  +2\nu\leq2$. We now use these facts to prove in the next lemma
estimates for $G$, which we shall invoke repeatedly in
the sequel to analyze various singular integrals. For the sake of clarity we
list those inequalities by their chronological order of appearance in the
proofs below.

\begin{lemma}
\label{l1}
Assume that Hypothesis ($K_{\beta,\beta^{\prime}}$) holds.
 Then, for all 
$x,y\in D$ and for every $\delta\in\left(\frac{d}{d+2},1\right)$ we have the
following inequalities.
\begin{description}
\item{(i)} For all $t, \tau, \sigma \in[0,T]$ with $t>\tau>\sigma$ and some $t^*\in(\sigma,\tau)$,
\begin{align}
&  \left|  G(x,t;y,\tau)-G(x,t;y,\sigma)\right| \nonumber\\
&  \leq c\left(  t-\tau\right)  ^{-\delta}(\tau-\sigma)^{\delta}(t-t^{\ast
})^{-\frac{d}{2}}\exp\left[  -c\frac{\left|  x-y\right|  _{{}}^{2}}{t-t^{\ast
}}\right].  \label{14}
\end{align}
\item{(ii)} For all $t,s,\tau\in[0,T]$ with $t>s>\tau$ and some $\tau^*\in(s,t)$,
\begin{align}
&  \left|  G(x,t;y,\tau)-G(x,s;y,\tau)\right|  ^{{}}\nonumber\\
&  \leq c\left(  t-s\right)  ^{\delta}(s-\tau)^{-\delta}(\tau^{\ast}
-\tau)^{-\frac{d}{2}}\exp\left[  -c\frac{\left|  x-y\right|  _{{}}^{2}}
{\tau^{\ast}-\tau}\right]  \label{15}
\end{align}
and
\begin{align}
&  \left|  G(x,t;y,\tau)-G(x,s;y,\tau)\right|  ^{\delta}\nonumber\\
&  \leq c\left(  t-s\right)  ^{\delta}(s-\tau)^{-\frac{d+2}{2}\delta+\frac
{d}{2}}(\tau^{\ast}-\tau)^{-\frac{d}{2}}\exp\left[  -c\frac{\left|
x-y\right|  _{{}}^{2}}{\tau^{\ast}-\tau}\right].  \label{16}
\end{align}
\item{(iii)} For all $t, s, \tau, \sigma\in[0,T]$ with $t>s>\tau>\sigma$, 
\begin{equation}
\left|  G(x,t;y,\tau)-G(x,t;y,\sigma)\right|  ^{1-\delta}
\leq c\left(  \tau-\sigma\right)  ^{1-\delta}(s-\tau)^{-\frac{d+2}
{2}(1-\delta)} \label{17}
\end{equation}
uniformly in $t$.
\end{description}
\end{lemma}

\noindent\textbf{Proof. }By applying successively (\ref{8}), the
mean-value theorem for $G$ and (\ref{13}) with $\left|
\mu\right|  =0$ and $\nu=1$ we may write
\begin{align*}
&  \left|  G(x,t;y,\tau)-G(x,t;y,\sigma)\right| \\
&  \leq\left(  \left|  G(x,t;y,\tau)\right|  +\left|  G(x,t;y,\sigma)\right|
\right)  ^{1-\delta}\left|  G(x,t;y,\tau)-G(x,t;y,\sigma)\right|  ^{\delta}\\
&  \leq c\left(  (t-\tau)^{-\frac{d}{2}}+(t-\sigma)^{-\frac{d}{2}}\right)
^{1-\delta}(\tau-\sigma)^{\delta}\left|  G_{t^{\ast}}(x,t;y,t^{\ast})\right|
^{\delta}\\
&  \leq c(t-\tau)^{-\frac{d}{2}(1-\delta)}(t-t^{\ast})^{-\frac{d+2}{2}
\delta+\frac{d}{2}}(\tau-\sigma)^{\delta}(t-t^{\ast})^{-\frac{d}{2}}
\exp\left[  -c\frac{\left|  x-y\right|  _{{}}^{2}}{t-t^{\ast}}\right] \\
&  \leq c\left(  t-\tau\right)  ^{-\delta}(\tau-\sigma)^{\delta}(t-t^{\ast
})^{-\frac{d}{2}}\exp\left[  -c\frac{\left|  x-y\right|  _{{}}^{2}}{t-t^{\ast
}}\right]
\end{align*}
for some $t^{\ast}\in(\sigma,\tau)$, since $-\frac{d+2}{2}\delta+\frac{d}{2}<0$
and $-\frac{d}{2}(1-\delta)-\frac{d+2}{2}\delta+\frac{d}{2}=-\delta$. This
proves (\ref{14}). Up to some minor but important changes, the
remaining inequalities can all be proved in a similar way. \hfill $\blacksquare$

\bigskip

Estimate (\ref{14}) now allows us to prove that our notion of
mild solution in Definition \ref{l2} is indeed well-defined; to this end for arbitrary mappings $\varphi$ and $u$
defined on $D$ and $D\times [0,T]$, respectively, we introduce
the functions $A(\varphi)$, $B(u)$, $C(u): D\times\left[0,T\right]  \mapsto\mathbb{R}$ by
\begin{align}
A(\varphi)(x,t)   &  :=\mathit{\ }\int_{D}dy\ G(x,t;y,0)\varphi
(y),\label{19}\\
B(u)(x,t)  &  :=\mathit{\ }\int_{0}^{t}d\tau\int_{D}dy\ G(x,t;y,\tau
)g\left(  u(y,\tau)\right)  ,\label{20}\\
C (u)(x,t)  &  :=\sum_{i=1}^{+\infty}\lambda_{i}^{\frac{1}{2}
}\int_{0}^{t}\left(  \int_{D}dy\ G(x,t;y,\tau)h\left(  u(y,\tau)\right)
e_{i}(y)\right)  B_{i}^{H}(d\tau), \label{21}
\end{align}
and prove the following result.

\begin{lemma}
\label{l2}
The hypotheses are the same as in Definition \ref{d2}.
 Then, for every
 $u\in\mathcal{B}^{\alpha,2}(0,T;L^{2}(D))$ we have $A(\varphi)(.,t) , B(u)(.,t)\in L^{2}(D)$, and also
$C(u)(.,t)\in L^{2}(D)$ a.s., for every $t\in\left[0,T\right]$.
\end{lemma}

\bigskip

\noindent\textbf{Proof. }The assertion is evident for  $A(\varphi)(.,t)$, 
since $\varphi$ is bounded and (\ref{8}) holds.
As for $B(u)(.,t)$, we infer from the Gaussian property of $G$ that
the measure $d\tau dy\left|  G(x,t;y,\tau)\right|  $ is finite on $\left[
0,T\right]  \times D$ uniformly in $(x,t)\in D\times\left[  0,T\right]  $, so
that by using successively Schwarz inequality with respect to this measure
along with Hypothesis (L) for $g$ we obtain
\begin{align*}
&  \left|B(u)(x,t)\right|  \leq\int_{0}^{t}d\tau\int_{D}dy\left|
G(x,t;y,\tau)g\left(  u(y,\tau)\right)  \right| \\
&  \leq c\left(  \int_{0}^{t}d\tau\int_{D}dy\left|  G(x,t;y,\tau)\right|
\left(  1+\left|  u(y,\tau)\right|  ^{2}\right)  \right)  ^{\frac{1}{2}}
\end{align*}
for every $x\in D$. We then get the inequalities
\begin{align*}
&  \left\|B(u)(.,t)\right\|  _{2}^{2}=\int_{D}dx\left|  \int
_{0}^{t}d\tau\int_{D}dy\ G(x,t;y,\tau)g\left(  u(y,\tau)\right)  \right|  ^{2}\\
&  \leq c\int_{0}^{t}d\tau\int_{D}dy\left(  1+\left|  u(y,\tau)\right|
^{2}\right) 
\leq c\left(  1+\int_{0}^{t}d\tau\left\|  u(.,\tau)\right\|  _{2}
^{2}\right) < +\infty. 
\end{align*}

 It remains to show that $\left\|C(u)(.,t)\right\|  _{2}^{2} <+\infty$ a.s. for every $t\in\left[  0,T\right] $. 

Define the functions $f_{i,t}(u):\left[0,t\right)  \mapsto L^{2}(D)$ by
\begin{equation}
f_{i,t}(u)(.,\tau):=\int_{D}dy\ G(.,t;y,\tau)h\left(  u(y,\tau)\right)
e_{i}(y). \label{22}
\end{equation}
We shall prove that 
\begin{equation}
\sum_{i=1}^{+\infty}\lambda_{i}^{\frac{1}{2}}\left\|  \int_{0}^{t}
f_{i,t}(u)(.,\tau)B_{i}^{H}(d\tau)\right\|  _{2}
 \leq c\ r_{\alpha}^{H}\left(  1+\left\|  u\right\|
_{\alpha,2,T}\right),  \label{23}
\end{equation}
a.s., where $r_\alpha^H$ is the a.s.  finite and positive random variable
defined in (\ref{29}).

Indeed, by using an argument similar to the one above, since $h$ is Lipschitz
continuous and $\sup_{i\in\mathbb{N}^{+}}\left\|  e_{i}\right\| _{\infty
}<+\infty$, we first obtain
\begin{equation}
\sup_{i\in\mathbb{N}^{+}}\left\|  f_{i,t}(u)(.,\tau)\right\|  _{2}\leq c\left(  1+\left\|
u(.,\tau)\right\|  _{2}\right)  \label{24}
\end{equation}
for every $\tau\in\left[0,t\right)$. Furthermore,
for every $x\in D$ and all $\sigma,\tau\in\left[  0,t\right)  $ with
$\tau>\sigma$ we have
\begin{align*}
&\left|  f_{i,t}(u)(x,\tau)-f_{i,t}(u)(x,\sigma)\right| 
 \leq c\Big(\int_{D}dy\left|  G(x,t;y,\tau)\right|  \left|  u(y,\tau
)-u(y,\sigma)\right| \\
&\quad \quad  +\int_{D}dy\left|  G(x,t;y,\tau)-G(x,t;y,\sigma)\right|  (1+\left|
u(y,\sigma)\right|)\Big),
\end{align*}
so that we get successively
\begin{align*}
&  \left|  f_{i,t}(u)(x,\tau)-f_{i,t}(u)(x,\sigma)\right|  ^{2}\\
&  \leq c\int_{D}dy\left|  G(x,t;y,\tau)\right|  \left|  u(y,\tau
)-u(y,\sigma)\right|  ^{2}\\
&  +c\int_{D}dy\left|  G(x,t;y,\tau)-G(x,t;y,\sigma)\right|  \left(  1+\left|
u(y,\sigma)\right|  ^{2}\right) \\
&  \leq c\int_{D}dy\left|  G(x,t;y,\tau)\right|  \left|  u(y,\tau
)-u(y,\sigma)\right|  ^{2}\\
&  +c\left(  t-\tau\right)  ^{-\delta}(\tau-\sigma)^{\delta}\int
_{D}dy(t-t^{\ast})^{-\frac{d}{2}}\exp\left[  -c\frac{\left|  x-y\right|  _{{}
}^{2}}{t-t^{\ast}}\right]  \left(  1+\left|  u(y,\sigma)\right|  ^{2}\right)
\end{align*}
for some $t^{\ast}\in(\sigma,\tau)$ and for every $\delta\in\left(  \frac
{d}{d+2},1\right)$. This is achieved by using Schwarz inequality with respect to the finite
measures $dy\left|  G(x,t;y,\tau)\right|  $ and $dy\left|  G(x,t;y,\tau
)-G(x,t;y,\sigma)\right|  $ on $D$, respectively, along with
(\ref{14}).  We then integrate the preceding estimate with
respect to $x\in D$ and apply the Gaussian property of $G$ to eventually
obtain
\begin{align}
&\sup_{i\in\mathbb{N}^{+}}\left\|  f_{i,t}(u)(.,\tau)-f_{i,t}(u)(.,\sigma)\right\|_{2}\nonumber\\
&\leq c\left(\left\|  u(.,\tau)-u(.,\sigma)\right\|  _{2}
+\left(  t-\tau\right)  ^{-\frac{\delta}{2}}(\tau-\sigma)^{\frac{\delta}{2}}\left(  1+\left\|  u(.,\sigma)\right\|  _{2}\right)\right).\label{25}
\end{align}
Therefore, by applying (\ref{i}) we have
\begin{align}
&  \sum_{i=1}^{+\infty}\lambda_{i}^{\frac{1}{2}}\left\|  \int_{0}^{t}
f_{i,t}(u)(.,\tau)B_{i}^{H}(d\tau)\right\|  _{2}\nonumber\\
&\le r_\alpha^H\sup_{i\in\mathbb{N}^+} \int_{0}^{t}d\tau\left(  \frac{\left\|
f_{i,t}(u)(.,\tau)\right\|  _{2}}{\tau^{\alpha}}+\int_{0}^{\tau}
d\sigma\ \frac{\left\|  f_{i,t}(u)(.,\tau)-f_{i,t}(u)(.,\sigma)\right\|  _{2}
}{\left(  \tau-\sigma\right)  ^{\alpha+1}}\right) \nonumber\\
&\le c r_\alpha^H\left(  1+\int_{0}^{t}d\tau\ \frac{\left\|
u(.,\tau)\right\|  _{2}^{{}}}{\tau^{\alpha}}+\int_{0}^{t}d\tau\int_{0}^{\tau
}d\sigma\ \frac{\left\|  u(.,\tau)-u(.,\sigma)\right\|  _{2}^{{}}}{\left(
\tau-\sigma\right)  ^{\alpha+1}}\right. \nonumber\\
&  +\left.  \int_{0}^{t}d\tau(t-\tau)^{-\frac{\delta}{2}}\int_{0}^{\tau
}d\sigma(\tau-\sigma)^{\frac{\delta}{2}-\alpha-1}\left(  1+\left\|
u(.,\sigma)\right\|  _{2}\right)  \right)  \label{26}
\end{align}
a.s.  
 
Let us now examine more closely the singular integrals in the
above terms. On the one hand, we may write
\begin{equation}
\int_{0}^{t}d\tau\ \frac{\left\|u(.,\tau)\right\|_2}{\tau^{\alpha}}
+\int_{0}^{t}d\tau\int_{0}^{\tau}d\sigma\ \frac{\left\|u(.,\tau)-u(.,\sigma)\right\|_2}{\left(\tau-\sigma\right)^{\alpha+1}}
\leq c\,\left\|  u\right\|_{\alpha,2,T},  \label{27}
\end{equation}
by using Schwarz inequality relative to the measure $d\tau$ on $(0,t)$ in
the last two integrals along with (\ref{3}). On the other hand, 
in (\ref{26}) the exponent $\delta$ can be taken arbitrarly close to $1$; consequently 
our range of values of $\alpha$ allows the condition $2\alpha<\delta$ to be satisfied.
Thus we can
integrate the singularities of the time increments in the last line of
(\ref{26}) and get the bound
\begin{equation}
\int_{0}^{t}d\tau(t-\tau)^{-\frac{\delta}{2}}\int_{0}^{\tau}d\sigma
(\tau-\sigma)^{\frac{\delta}{2}-\alpha-1}\left(  1+\left\|  u(.,\sigma)\right\|  _{2}\right)
\leq c\left(1+\sup_{t\in[0,T]}\left\|u(.,t)\right\|_2\right).
\label{28}
\end{equation}
Therefore, we can
substitute (\ref{27}), (\ref{28}) into (\ref{26}) to
obtain (\ref{23}). 

\hfill $\blacksquare$

\bigskip

In order to relate the notions of variational and mild solution, we recall
that the self-adjoint operator $A(t)=-\operatorname*{div}(k(.,t)\nabla)$
defined on (\ref{11}) generates the family of evolution operators
$U(t,s)_{0\leq s\leq t\leq T}$ in $L^{2}(D)$ given by
\begin{equation}
U(t,s)v=\begin{cases}
v,&\text{if}\ \ s=t,\\
\int_{D}dy\ G(.,t;y,s)v(y), &\text{if}\ \ t>s,
\label{30}
\end{cases}
\end{equation}
and that each such $U(t,s)$ is itself self-adjoint (see, for instance,
\cite{tanabe}), which means that the symmetry property 
\begin{equation}
G(x,t;y,s)=G(y,t;x,s) \label{symmetry}
\end{equation}
holds for every $(x,t;y,s)\in\overline{D}\times\left[  0,T\right]
\times\overline{D}\times\left[  0,T\right]  \diagdown\left\{  s,t\in\left[0,T\right]  :s\geq t\right\}$. 

\bigskip

\noindent{\bf Proof of Statement (a) of Theorem \ref{t1}}
\medskip

The existence of a variational solution $u_V$ was proved in Theorem of \cite{nuavu}.

In order to prove that every variational solution is mild, we follow the same approach as in Theorem 2 of \cite{sansovu1}.
For the sake of completeness, we sketch the main ideas. 

We shall check that the $L^{2}(D)$-valued stochastic process
\begin{align*}
&  u_{V}(.,t)-\int_{D}dy\ G(.,t;y,0)\varphi(y)-\int_{0}^{t}d\tau\int
_{D}dy\ G(.,t;y,\tau)g\left(  u_{V}(y,\tau)\right) \\
&  -\sum_{i=1}^{+\infty} \lambda_i^{\frac{1}{2}}\int_{0}^{t}\left(\int_D dy\ G(.,t;y,\tau)h\left(u_{V}(y,\tau)e_i(y)\right)\right)
B_i^{H}(d\tau)
\end{align*}
is a.s. orthogonal for every $t\in\left[  0,T\right]  $ to the dense subspace
$\mathcal{C}_{0}^{2}(D)$ consisting of all twice continuously differentiable
functions with compact support in $D$. To this end, for every $v\in\mathcal{C}_{0}^{2}(D)$
and all $s,t\in\left[  0,T\right]$ with $t\geq s$
we define $v^{t}(.,s):=U(t,s)v$, that is, 
\begin{equation}
v^{t}(x,s)=
\begin{cases}
v(x),& \text{if}\ \ s=t,\\
\int_{D}dyG(y,t;x,s)v(y),& \text{if}\ \ t>s,
\end{cases}
\label{31}
\end{equation}
for every $x\in D$ by taking (\ref{30}) and (\ref{symmetry}) into
account. It then follows from (\ref{12}), (\ref{symmetry}) and Gauss'
divergence theorem that $v^{t}\in H^{1}(D\times(0,T))$, and that for every $t\in\left[  0,T\right]$, the relation
\begin{equation}
\int_{0}^{t}d\tau\int_{D}dx\ v_{\tau}^{t}(x,\tau)u_{V}(x,\tau)=\int_{0}^{t}
d\tau\int_{D}dx\left(  \nabla v^{t}(x,\tau),k(x,\tau)\nabla u_{V}
(x,\tau)\right)  _{\mathbb{R}^{d}} \label{32}
\end{equation}
holds a.s.
Therefore, we may take (\ref{31}) as a test function in
(\ref{4}), which, as a consequence of (\ref{32}),
leads to the relation
\begin{align*}
(v,u_{V}(.,t))_{2}  &  =(v^{t}(.,0),\varphi)_{2}+\int_{0}^{t}d\tau
(v^{t}(.,\tau),g(u_{V}(.,\tau)))_{2}\\
&  +\sum_{i=1}^{+\infty}\lambda_{i}^{\frac{1}{2}}\int_{0}^{t}\left(
v^{t}(.,\tau),h(u_{V}(.,\tau))e_{i}\right)  _{2}B_{i}^{\text{\textsc{H}}%
}(d\tau),
\end{align*}
valid a.s. for every $t\in\left[0,T\right]$. After some rearrangements, the substitution of
(\ref{31}) into the right-hand side of the preceding expression then
leads to the equality
\begin{align*}
(v,u_{V}(.,t))_{2}  &  =\left(  v,\int_{D}dyG(.,t;y,0)\varphi(y)\right)_{2}\\
&+\left(  v,\int_{0}^{t}d\tau\int_{D}dyG(.,t;y,\tau)g\left(  u_{V}
(y,\tau)\right)  \right)_{2}\\
&  +\left(  v,\sum_{i=1}^{+\infty} \lambda_i^{\frac{1}{2}}\int_{0}^{t}\left(\int_D dy\  G(.,t;y,\tau)h\left(  u_{V}(y,\tau)e_i(y)\right)\right)
B_i^{H}(d\tau)\right)_{2},
\end{align*}
which holds for every $t\in\left[  0,T\right]$ a.s.  and every
$v\in\mathcal{C}_{0}^{2}(D)$, thereby leading to the desired orthogonality
property. \hfill $\blacksquare$
\bigskip
\newpage

\noindent{\bf Proof of Statement $(b)$ of Theorem \ref{t1}}
\bigskip

\noindent Under the standing assumptions, we already know from \cite{nuavu}
that the variational solution is unique. Moreover, we have just proved that every variational solution is also
a mild solution. Hence, it suffices to prove that uniqueness holds
within the class of mild solutions.
To this
end, let us write $u_{M}$ and $\tilde{u}_{M}$ for any two such solutions
corresponding to the same initial condition $\varphi$; from (\ref{9}) and
(\ref{19})-(\ref{21}) we have
\begin{align}
&  \left\Vert u_{M}(.,t)-\tilde{u}_{M}(.,t)\right\Vert _{2}\nonumber\\
&  \leq\left\Vert B(u_{M})(.,t)-B(\tilde{u}_{M})(.,t)\right\Vert
_{2}+\left\Vert C(u_{M})(.,t)-C(\tilde{u}_{M})(.,t)\right\Vert _{2}
\label{302}
\end{align}
a.s. for every $t\in\left[  0,T\right]  $.

We proceed by estimating both terms
on the right-hand side of (\ref{302}). 
Since $g$ is Lipschitz, we have
\begin{equation}
\left\| B\left(  u_M\right)  (.,t)-B(\tilde{u}_M)(.,t)\right\|  _{2}^{2}
\leq c\int_{0}^{t}d\tau\left\|  u_M(.,\tau)-\tilde{u}_M(.,\tau)\right\|  _{2}^{2} \label{303}
\end{equation}
a.s. for every $t\in\left[ 0,T\right]$.

In order to analyze the second term or the right-hand side of (\ref{302}) we will need the following preliminary result.

\begin{lemma}
\label{l5}
The hypotheses are the same as in part (a) of Theorem \ref{t1} and let
the $f_{i,t}(u)$'s be the functions given by (\ref{22}). Then,
the estimate
\begin{equation}
\sup_{(i,t)\in\mathbb{N}^+\times[0,T]} \left\| f_{i,t}(u_M)(.,\tau)-f_{i,t}(\tilde{u}_M)(.,\tau)\right\|_2
\leq c\left\|  u_M(.,\tau)-\tilde{u}_M(.,\tau)\right\|_{2}
\label{304}
\end{equation}
holds a.s. for every $\tau\in\left[0,t\right)$.

Moreover, if $h$ is an affine function we
have
\begin{align}
&  \sup_{i\in\mathbb{N}^+}\left\|  \ f_{i,t}(u_M)(.,\tau)-f_{i,t}(\tilde{u}_M)(.,\tau)
 -f_{i,t}(u_M)(.,\sigma)+f_{i,t}(\tilde{u}_M)(.,\sigma)\right\|_{2}\nonumber\\
&\quad  \leq c(t-\tau)^{-\delta}(\tau-\sigma)^{\delta}\left\|
u_M(.,\sigma)-\tilde{u}_M(.,\sigma)\right\|  _{2}\nonumber\\
& \quad +c\left\|  u_M(.,\tau)-\tilde{u}_M(.,\tau)-u_M(.,\sigma)
+\tilde{u}_M(.,\sigma)\right\|  _{2} \label{305}
\end{align}
a.s. for all $t,\tau,\sigma\in\left[0,T\right]$ with $t>\tau>\sigma$ and every $\delta\in\left(\frac{d}{d+2},1\right)$.
\end{lemma}

\noindent\textbf{Proof. }Up to minor modifications, we can prove (\ref{304}) as we
argued in the proof of (\ref{24}).

For the proof of (\ref{305}) we first write
\begin{align*}
&\left\|  \ f_{i,t}(u_M)(.,\tau)-f_{i,t}(\tilde{u}_M)(.,\tau)
 -f_{i,t}(u_M)(.,\sigma)+f_{i,t}(\tilde{u}_M)(.,\sigma)\right\|_{2}^2\\
 &\quad \le 2\left(F^1(i,t,\tau,\sigma)+F^2(i,t,\tau,\sigma)\right),
 \end{align*}
 with
 \begin{align*}
& F^1(i,t,\tau,\sigma)\\
&= \int_D dx\left\vert\int_D dy\ e_i(y)\left(G(x,t;y,\tau)-(G(x,t;y,\sigma)\right)\left(u_M(y,\tau)-\tilde u_M(y,\tau)\right)\right\vert^2,\\
&F^2(i,t,\tau,\sigma)\\
&=\int_D dx\left\vert\int_D dy\ e_i(y) G(x,t;y,\sigma)\left(u_M(y,\tau)-\tilde u_M(y,\tau)-u_M(y,\sigma)+\tilde u(y,\sigma)\right)\right\vert^2.
\end{align*}
From the Gaussian property of $G$ we clearly see that  $F^2(i,t,\tau,\sigma)$ is  bounded above by the square of the last term of (\ref{305}). 
Moreover, by applying first (\ref{14}) and then Schwarz inequality we obtain
\beqn
F^1(i,t,\tau,\sigma)\le c(t-\tau)^{-2\delta}(\tau-\sigma)^{2\delta}\left\|u_M(.,\sigma)-\tilde{u}_M(.,\sigma)\right\|_{2}^2.
\eeqn
Hence (\ref{305}) is proved.

\hfill $\blacksquare$
\medskip

The preceding result now leads to the following estimate for the second term
on the right-hand side of (\ref{302}).

\begin{lemma}
\label{l6}
The hypotheses are those of Theorem \ref{t1} part (b).
Then we have
\begin{align}
&  \left\|C(u_M)(.,t)-C(\tilde{u}_M
)(.,t)\right\|  _{2}\nonumber\\
&\quad  \leq cr_{\alpha}^{H}\left(  \int_{0}^{t}d\tau\left(
\frac{1}{\tau^{\alpha}}+\frac{1}{\left(  t-\tau\right)  ^{\alpha}}\right)
\left\|  u_M(.,\tau)-\tilde{u}_M(.,\tau)\right\|  _{2}\right. \nonumber\\
&\quad  +\left.  \int_{0}^{t}d\tau\int_{0}^{\tau}d\sigma\frac{\left\|  u_M
(.,\tau)-\tilde{u}_M(.,\tau)-u_M(.,\sigma)+\tilde{u}_M(.,\sigma
)\right\|  _{2}}{\left(  \tau-\sigma\right)  ^{\alpha+1}}\right)
\label{306}
\end{align}
a.s. for every $t\in\left[ 0,T\right]$.
\end{lemma}

\noindent\textbf{Proof. }From (\ref{21}), (\ref{22}), 
and by using (\ref{i}), (\ref{304}), (\ref{305}), we have
\begin{align}
&  \left\| C(u_M)(.,t)-C(\tilde{u}_M
)(.,t)\right\|  _{2}\nonumber\\
&  \leq cr_{\alpha}^{H}\left(  \int_{0}^{t}d\tau\frac{\left\|
u_M(.,\tau)-\tilde{u}_M(.,\tau)\right\|  _{2}}{\tau^{\alpha}}\right.
\nonumber\\
&  +\int_{0}^{t}d\tau\int_{0}^{\tau}d\sigma(t-\tau)^{-\delta}
(\tau-\sigma)^{\delta-\alpha-1}\left\|  u_M(.,\sigma
)-\tilde{u}_M(.,\sigma)\right\|  _{2}\nonumber\\
&  +\left.  \int_{0}^{t}d\tau\int_{0}^{\tau}d\sigma\frac{\left\|  u_M
(.,\tau)-\tilde{u}_M(.,\tau)-u_M(.,\sigma)+\tilde{u}_M(.,\sigma
)\right\|  _{2}}{\left(  \tau-\sigma\right)  ^{\alpha+1}}\right)
\label{307}
\end{align}
a.s. for every $t\in\left[  0,T\right]$.

Furthermore, by swapping each
integration variable for the other in the second term on the right-hand side
and by using Fubini's theorem we may write
\begin{align*}
&  \int_{0}^{t}d\tau\int_{0}^{\tau}d\sigma(t-\tau)^{-\delta}
(\tau-\sigma)^{\delta-\alpha-1}\left\|  u_M(.,\sigma
)-\tilde{u}_M(.,\sigma)\right\|  _{2}\\
&  =\int_{0}^{t}d\tau\left\|  u_M(.,\tau)-\tilde{u}_M(.,\tau)\right\|
_{2}\int_{\tau}^{t}d\sigma(t-\sigma)^{-\delta}(\sigma-\tau
)^{\delta-\alpha-1}\\
&  =c\int_{0}^{t}d\tau\frac{\left\|  u_M(.,\tau)-\tilde{u}_M
(.,\tau)\right\|  _{2}}{\left(  t-\tau\right)  ^{\alpha}},
\end{align*}
after having evaluated the singular integral explicitly in terms of Euler's
Beta function; this is possible since we can choose $\delta\in\left(\frac{d}{d+2},1\right)$ such that $\alpha<\delta$. The
substitution of the preceding expression into (\ref{307}) then proves
(\ref{306}). \hfill $\blacksquare$
\bigskip

In what follows, we write $R$ for all the  irrelevant a.s. finite and positive
random variables that appear in the different estimates, unless we specify these
variables otherwise.
The preceding inequalities then lead to a
crucial estimate for $z_{M}:=u_{M}-\tilde{u}_{M}$ with respect to the norm in
$B^{\alpha,2}(0,t;L^{2}(D))$.
\begin{lemma}
\label{l30}
We assume the same hypotheses as in part (b) of Theorem \ref{t1}.
Then we have
\begin{align*}
\left\Vert z_{M}\right\Vert _2^{2}  &  \leq R\left(  \int_{0}^{t}d\tau
\sup_{\sigma\in\lbrack0,\tau]}\left\Vert z_{M}(.,\sigma)\right\Vert_{2}^{2}\right. \nonumber\\
&  +\left.  \int_{0}^{t}d\tau\left(  \int_{0}^{\tau}d\sigma\frac{\left\Vert
z_{M}(.,\tau)-z_{M}(.,\sigma)\right\Vert _{2}}{\left(  \tau-\sigma\right)^{\alpha+1}}\right)^{2}\right)  
\end{align*}
a.s. for every $t\in\left[  0,T\right]$.
\end{lemma}

\noindent{\bf Proof.}
We apply Schwarz inequality relative to the measure $d\tau$ on
$(0,t)$ to both integrals on the right-hand side of (\ref{306}). This leads
to
\begin{align}
&  \Vert C(u_{M})(.,t)-C(\tilde{u}_{M})(.,t)\Vert_{2}^{2}\nonumber\\
&  \leq R\left(  \int_{0}^{t}d\tau\Vert z_{M}(.,\tau)\Vert_{2}^{2}+\int
_{0}^{t}d\tau\left(  \int_{0}^{\tau}d\sigma\frac{\Vert z_{M}(.,\tau
)-z_{M}(.,\sigma)\Vert_{2}}{(\tau-\sigma)^{\alpha+1}}\right)  ^{2}\right)
\label{3701}
\end{align}
a.s. for every $t\in\left[  0,T\right]  $. This estimate along with (\ref{302}),
(\ref{303}) yields the result.

 \hfill $\blacksquare$

As a consequence of the preceding Lemma we obtain
\begin{align}
\left\Vert z_{M}\right\Vert _{\alpha,2,t}^{2}  &  \leq R\left(  \int_{0}^{t}d\tau
\sup_{\sigma\in\lbrack0,\tau]}\left\Vert z_{M}(.,\sigma)\right\Vert_{2}^{2}\right. \nonumber\\
&  +\left.  \int_{0}^{t}d\tau\left(  \int_{0}^{\tau}d\sigma\frac{\left\Vert
z_{M}(.,\tau)-z_{M}(.,\sigma)\right\Vert _{2}}{\left(  \tau-\sigma\right)^{\alpha+1}}\right)^{2}\right)  \label{308},
\end{align}
by the very definition of the norm $\left\Vert \cdot\right\Vert _{\alpha,2,t}^{2}$.  

\bigskip

We proceed by analyzing further the second
term on the right-hand side of (\ref{308}), so as to eventually obtain an inequality of Gronwall type
for $\Vert z_M\Vert_{\alpha,2,t}^2$.  

First we introduce some notation.
For $0\le \tau<s\le t\le T$, we set
\begin{equation}
f_{i,t,s}^{\ast}(u_{M})(.,\tau):=f_{i,t}(u_{M})(.,\tau)-f_{i,s}(u_{M})(.,\tau),
\label{39}
\end{equation}
where the $f_{i,t}(u_{M})$'s are given by (\ref{22}).

By reference to (\ref{9}),
we may write
\begin{align}
&  z_M(.,\tau)-z_M(.,\sigma)\nonumber\\
&  =\int_{\sigma}^{\tau}d\rho\int_{D}dyG(.,\tau;y,\rho)\left(  g(u_M
(y,\rho))-g(\tilde{u}_M(y,\rho))\right) \nonumber\\
&  +\int_{0}^{\sigma}d\rho\int_{D}dy\left(  G(.,\tau;y,\rho)-G(.,\sigma
;y,\rho)\right)  \left(  g(u_M(y,\rho))-g(\tilde{u}_M
(y,\rho))\right) \nonumber\\
&  +\sum_{i=1}^{+\infty}\lambda_{i}^{\frac{1}{2}}\int_{\sigma}^{\tau}\left(
f_{i,\tau}(u_M)(.,\rho)-f_{i,\tau}(\tilde{u}_M
)(.,\rho)\right)  B_{i}^{H}(d\rho)\nonumber\\
&  +\sum_{i=1}^{+\infty}\lambda_{i}^{\frac{1}{2}}\int_{0}^{\sigma}\left(
f_{i,\tau,\sigma}^{\ast}(u_M)(.,\rho)-f_{i,\tau,\sigma}^{\ast
}(\tilde{u}_M)(.,\rho)\right)  B_{i}^{H}(d\rho)
\label{309}
\end{align}
for all $\sigma,\tau\in\left[  0,t\right]  $ with $\tau>\sigma$.

Our next goal is to estimate the $L^{2}(D)$-norm of each contribution on
the right-hand side of (\ref{309}). Regarding the first two terms
we have the following result. 

\begin{lemma}
\label{l7}
The hypotheses are the same as in part (a) of Theorem \ref{t1}; then we
have
\begin{align}
&  \left\|  \int_{\sigma}^{\tau}d\rho\int_{D}dy\ G(.,\tau;y,\rho)\left(
g(u_M(y,\rho))-g(\tilde{u}_M(y,\rho))\right)  \right\|
_{2}\nonumber\\
&  \leq c\left(  \tau-\sigma\right)  ^{\frac{1}{2}}\left(  \int_{\sigma}
^{\tau}d\rho\left\|  z_M(.,\rho)\right\|  _{2}^{2}\right)
^{\frac{1}{2}}
\label{310}
\end{align}
and
\begin{align}
&  \left\|  \int_{0}^{\sigma}d\rho\int_{D}dy\left(  G(.,\tau;y,\rho
)-G(.,\sigma;y,\rho)\right)  \left(  g(u_M(y,\rho
))-g(\tilde{u}_M(y,\rho))\right)  \right\|  _{2}\nonumber\\
&  \leq c\left(  \tau-\sigma\right)  ^{\frac{\delta}{2}}\left(  \int
_{0}^{\sigma}d\rho\left(  \sigma-\rho\right)  ^{-\delta}\left\|  z_M
(.,\rho)\right\|  _{2}^{2}\right)  ^{\frac{1}{2}}
\label{311}
\end{align}
a.s. for all $\sigma,\tau\in\left[  0,t\right] $ with 
$\tau>\sigma$ and every $\delta\in\left(  \frac{d}{d+2},1\right)$.
\end{lemma}
{\bf Proof:} The inequality (\ref{310}) follows by applying Schwarz inequality and using the Gaussian property along with
assumption (L). As for (\ref{311}), we first apply Schwarz inequality  with respect to the measure on $[0,\sigma]\times D$ given by
$|G(x,\tau;y,\rho)-G(x,\sigma;y,\rho)|d\rho\ dy$ and then (\ref{15}).

\hfill $\blacksquare$

Next, we turn to the analysis of the third term on the right-hand side of
(\ref{309}).

\begin{lemma}
\label{l8}
With the same hypotheses as in part (b) of Theorem \ref{t1}, we have
\begin{align*}
&  \sum_{i=1}^{+\infty}\lambda_{i}^{\frac{1}{2}}\left\|  \int_{\sigma}^{\tau
}\left(  f_{i,\tau}(u_M)(.,\rho)-f_{i,\tau}(\tilde{u}_M
)(.,\rho)\right)  B_{i}^{H}(d\rho)\right\|  _{2}\\
&  \leq R\left(  \int_{\sigma}^{\tau}
d\rho\left(  \frac{1}{\left(  \rho-\sigma\right)  ^{\alpha}}+\frac{1}{\left(
\tau-\rho\right)  ^{\alpha}}\right)  \left\|  z_M(.,\rho)\right\|
_{2}\right.  \\
&  +\left.  \int_{\sigma}^{\tau}d\rho\int_{\sigma}^{\rho}d\varsigma
\frac{\left\|  z_M(.,\rho)-z_M(.,\varsigma)\right\|  _{2}
}{\left(  \rho-\varsigma\right)  ^{\alpha+1}}\right)
\end{align*}
a.s. for all $\sigma,\tau\in\left[0,t\right]$ with 
$\tau>\sigma$.
\end{lemma}
\noindent\textbf{Proof.} 
In terms of the variables $\tau,\rho$ and $\varsigma$,
inequalities (\ref{304}), (\ref{305}) of Lemma \ref{l5} now read
\begin{equation}
\sup_{(i,\tau)\in\mathbb{N}^+\times [0,T]}\left\|  \ f_{i,\tau}(u_M)(.,\rho)-f_{i,\tau}
(\tilde{u}_M)(.,\rho)\right\|  _{2}
\leq c\left\|  z_M(.,\rho)\right\|_{2} \label{312}
\end{equation}
and
\begin{align}
& \sup_{i\in\mathbb{N}^+} \left\|  \ f_{i,\tau}(u_M)(.,\rho)-f_{i,\tau}
(\tilde{u}_M)(.,\rho)\ -f_{i,\tau}(u_M)(.,\varsigma
)+f_{i,\tau}(\tilde{u}_M)(.,\varsigma)\right\|  _{2}^{{}}\nonumber\\
&  \leq c(\tau-\rho)^{-\delta}(\rho-\varsigma)^{\delta}
\left\|  z_M(.,\varsigma)\right\|  _{2}
+c\left\|  z_M(.,\rho)-z_M(.,\varsigma)\right\|  _{2},
\label{313}
\end{align}
respectively. Thus, by an extended version of (\ref{i}) for indefinite generalized Stieltjes integrals
(see Proposition 4.1 in \cite{nuarasca}),
\begin{align*}
&  \sum_{i=1}^{+\infty}\lambda_{i}^{\frac{1}{2}}\left\|  \int_{\sigma}^{\tau}
\left(  f_{i,\tau}(u_M)(.,\rho)-f_{i,\tau}(\tilde{u}_M)(.,\rho)\right)  B_{i}^{H}(d\rho)\right\|  _{2}\\
&  \leq r_\alpha^H \sup_{i\in\mathbb{N}^+}\Big(
\Big(\int_{\sigma}^{\tau}d\rho\frac{\left\|\ f_{i,\tau}(u_M)(.,\rho)-f_{i,\tau}(\tilde{u}_M%
)(.,\rho)\right\|_{2}}{\left(  \rho-\sigma\right)  ^{\alpha}}\\
&  +\int_{\sigma}^{\tau}d\rho\int_{\sigma}^{\rho}\frac{d\varsigma}{\left(\rho-\varsigma\right)^{\alpha+1}}\\
&\times\left\|  \ f_{i,\tau}(u_M)(.,\rho)-f_{i,\tau}(\tilde{u}_M)(.,\rho)\ -f_{i,\tau}(u_M)(.,\varsigma)+f_{i,\tau}(\tilde{u}_M)(.,\varsigma)\right\|_2\Big)\Big) \\
&  \leq R\left(  \int_{\sigma}^{\tau}d\rho
\frac{\left\|  z_M(.,\rho)\right\|  _{2}}{\left(  \rho
-\sigma\right)  ^{\alpha}}+\int_{\sigma}^{\tau}d\rho\left(  \tau-\rho\right)
^{-\delta}\int_{\sigma}^{\rho}d\varsigma\left(  \rho-\varsigma
\right) ^{\delta-\alpha-1}\left\|  z_M(.,\varsigma
)\right\|  _{2}\right. \\
&  +\left.  \int_{\sigma}^{\tau}d\rho\int_{\sigma}^{\rho}d\varsigma
\frac{\left\|  z_M(.,\rho)-z_M(.,\varsigma)\right\|_{2}
}{\left(  \rho-\varsigma\right)  ^{\alpha+1}}\right)
\end{align*}
a.s. for all $\sigma,\tau\in\left[  0,t\right]  $ with $\tau>\sigma$ and every
$\delta\in\left(  \frac{d}{d+2},1\right)$. But the second term on the
right-hand side is equal to
\begin{equation*}
c\int_{\sigma}^{\tau}d\rho\left(  \tau-\rho\right)  ^{-\alpha}\left\|
z_M(.,\rho)\right\|  _{2},
\end{equation*}
as can be easily checked by applying Fubini's theorem and by evaluating the
resulting inner integral in terms of Euler's Beta function. \hfill $\blacksquare$
\bigskip

As for the analysis of the fourth term on the right-hand side of
(\ref{309}) we need some preparatory results. In particular we shall use the estimate for 
time increments of the Green function, valid for any  $\delta\in\left(\frac{d}{d+2},1\right)$: 
\begin{align}
& \left|G(x,t;y,\tau)-G(x,s;y,\tau)-G(x,t;y,\sigma)+G(x,s;y,\sigma)\right|\nonumber\\
& \leq \left( \left| G(x,t;y,\tau)-G(x,s;y,\tau)\right|^{\delta}
+\left|G(x,t;y,\sigma)-G(x,s;y,\sigma)\right|^{\delta}\right) \nonumber\\
&  \times\left( \left| G(x,t;y,\tau)-G(x,t;y,\sigma)\right|^{1-\delta}
+\left|G(x,s;y,\tau)-G(x,s;y,\sigma)\right|^{1-\delta}\right) \nonumber\\
&  \leq \left( t-s\right)^\delta(s-\tau)^{-\frac{d+2}{2}\delta+\frac{d}{2}}
\left( \tau-\sigma\right) ^{1-\delta}(s-\tau)^{-\frac{d+2}{2}(1-\delta)}\nonumber\\
& \times\left( (\tau^{\ast}-\tau)^{-\frac{d}{2}}
\exp\left[ -c\frac{\left|x-y\right|^{2}}{\tau^{\ast}-\tau}\right]  
+(\sigma^{\ast}-\sigma)^{-\frac{d}{2}}\exp\left[-c\frac{\left| x-y\right|^{2}}{\sigma^{\ast}-\sigma}\right]  \right) \nonumber\\
&  =c\left(  t-s\right)  ^{\delta}\left(  s-\tau\right)  ^{-1}\left(
\tau-\sigma\right)  ^{1-\delta}\nonumber\\
&  \times\left(  (\tau^{\ast}-\tau)^{-\frac{d}{2}}\exp\left[  -c\frac{\left|
x-y\right|  _{{}}^{2}}{\tau^{\ast}-\tau}\right]  +(\sigma^{\ast}
-\sigma)^{-\frac{d}{2}}\exp\left[  -c\frac{\left|  x-y\right|^{2}}{\sigma^{\ast}-\sigma}\right]  \right),  \label{44}
\end{align}
with $\tau^*, \sigma^* \in(s,t)$, which follows from (\ref{16})--(\ref{17}).

\begin{lemma}
\label{l9}
The hypotheses are the same as in part (b) of Theorem \ref{t1} and the
$f_{i,\tau,\sigma}^{\ast}(u)$'s are the functions given by
(\ref{39}). Then, the estimates
\begin{equation}
\sup_{i\in\mathbb{N}^+}\left\|\mathit{\ }f_{i,\tau,\sigma}^{\ast}(u_M)
(.,\rho)-\mathit{\ }f_{i,\tau,\sigma}^{\ast}(\tilde{u}_M
)(.,\rho)\right\|  _{2}
\leq c(\tau-\sigma)^{\frac{\delta}{2}}(\sigma-\rho)^{-\frac{\delta}{2}}\left\|z_M(.,\rho)\right\|  _{2} \label{314}
\end{equation}
and
\begin{align}
& \sup_{i\in\mathbb{N}^+}\left\|\mathit{\ }f_{i,\tau,\sigma}^{\ast}(u_M
)(.,\rho)-\mathit{\ }f_{i,\tau,\sigma}^{\ast}(\tilde{u}_M
)(.,\rho)-\mathit{\ }f_{i,\tau,\sigma}^{\ast}(u_M)(.,\varsigma
)+\mathit{\ }f_{i,\tau,\sigma}^{\ast}(\tilde{u}_M)(.,\varsigma
)\right\|  _{2}\nonumber\\
&  \leq c(\tau-\sigma)^{\frac{\delta}{2}}\left(  (\sigma-\rho)^{-\frac{1}{2}
}(\rho-\varsigma)^{\frac{1}{2}(1-\delta)}\left\|  z_M(.,\varsigma
)\right\|  _{2}\right. \nonumber\\
&  +\left.  (\sigma-\rho)^{-\frac{\delta}{2}}\left\|  z_M
(.,\rho)-z_M(.,\varsigma)\right\|  _{2}\right)  \label{315}
\end{align}
hold a.s. for all $\tau,\sigma,\rho,\varsigma\in\left[ 0,T\right]$ with 
$\tau>\sigma>\rho>\varsigma$ and every $\delta\in\left(  \frac{d}{d+2},1\right)$.
\end{lemma}
\noindent\textbf{Proof.} It follows from the same type of 
arguments as those outlined in the
proof of Lemma \ref{l5}. For the proof of (\ref{314}) the key estimate is (\ref{15}). For (\ref{315}), we also apply (\ref{15}) along with
(\ref{44}). 

 \hfill $\blacksquare$
\bigskip

The last relevant $L^{2}(D)$-estimate regarding (\ref{309}) is then
the following.
\begin{lemma}
\label{l10}
The hypotheses are the same as in part (b) of Theorem \ref{t1}. Then we
have
\begin{align*}
&  \sum_{i=1}^{+\infty}\lambda_{i}^{\frac{1}{2}}\left\|  \int_{0}^{\sigma
}\left(  f_{i,\tau,\sigma}^{\ast}(u_M)(.,\rho)-f_{i,\tau,\sigma
}^{\ast}(\tilde{u}_M)(.,\rho)\right)  B_{i}^{\text{\textsc{H}}
}(d\rho)\right\|  _{2}\\
&  \leq R (\tau-\sigma)^{\frac{\delta}{2}}\left(
\int_{0}^{\sigma}d\rho(\sigma-\rho)^{-\frac{\delta}{2}}\left(  \frac{1}
{\rho^{\alpha}}+\frac{1}{(\sigma-\rho)^{\alpha}}\right)  \left\|
z_M(.,\rho)\right\|  _{2}\right. \\
&  +\left.  \int_{0}^{\sigma}d\rho(\sigma-\rho)^{-\frac{\delta}{2}}\int
_{0}^{\rho}d\varsigma\frac{\left\|  z_M(.,\rho)-z_M
(.,\varsigma)\right\|  _{2}}{\left(  \rho-\varsigma\right)  ^{\alpha+1}
}\right)
\end{align*}
a.s. for all $\sigma,\tau\in\left[  0,t\right]$ with 
$\tau>\sigma$ and every $\delta\in\left(  \frac{d}{d+2},1-2\alpha\right)$.
\end{lemma}
\noindent\textbf{Proof.} By applying (\ref{i}), 
together with (\ref{314}), (\ref{315}), we get
\begin{align*}
&  \sum_{i=1}^{+\infty}\lambda_{i}^{\frac{1}{2}}\left\|  \int_{0}^{\sigma
}\left(  f_{i,\tau,\sigma}^{\ast}(u_M)(.,\rho)-f_{i,\tau,\sigma
}^{\ast}(\tilde{u}_M)(.,\rho)\right)  B_{i}^{\text{\textsc{H}}
}(d\rho)\right\|  _{2}\\
&  \leq R (\tau-\sigma)^{\frac{\delta}{2}}\left(
\int_{0}^{\sigma}d\rho(\sigma-\rho)^{-\frac{\delta}{2}}\frac{\left\|
z_M(.,\rho)\right\|  _{2}}{\rho^{\alpha}}\right. \\
&  +\int_{0}^{\sigma}d\rho(\sigma-\rho)^{-\frac{1}{2}}\int_{0}^{\rho
}d\varsigma(\rho-\varsigma)^{\frac{1}{2}(1-\delta)-\alpha-1}\left\|
z_M(.,\varsigma)\right\|  _{2}\\
&  +\left.  \int_{0}^{\sigma}d\rho(\sigma-\rho)^{-\frac{\delta}{2}}\int
_{0}^{\rho}d\varsigma\frac{\left\|  z_M(.,\rho)-z_M(.,\varsigma)
\right\|  _{2}}{\left(  \rho-\varsigma\right)  ^{\alpha+1}
}\right)
\end{align*}
a.s. for all $\sigma,\tau\in\left[  0,t\right]  $ with $\tau>\sigma$. But for every 
$\delta\in\left(  \frac{d}{d+2},1-2\alpha\right)  $, we have 
\begin{align*}
& \int_{0}^{\sigma}d\rho(\sigma-\rho)^{-\frac{1}{2}}\int_{0}^{\rho}
d \varsigma(\rho-\varsigma)^{\frac{1}{2}(1-\delta)-\alpha-1}\left\|
z_M(.,\varsigma)\right\|  _{2}\\
&  =c\int_{0}^{\sigma}d\rho(\sigma-\rho)^{-\alpha-\frac{\delta}{2}}\left\|
z_M(.,\rho)\right\| _{2}.
\end{align*}
This yields the result.

\hfill$\blacksquare$
\bigskip

Let us go back to the inequality (\ref{308}). Owing to (\ref{309}) and by using the estimates
(\ref{310}), (\ref{311}) together with Lemmas \ref{l8} and \ref{l10}, we have
\begin{equation}
\label{3160}
\Vert z_M\Vert_{\alpha,2,t}^2 \le R \int_0^t d\tau \left(\sup_{\rho\in[0,\tau]}\Vert z_M(.,\rho)\Vert_2^2
+\sum_{k=1}^6 \left[I_k(\tau)\right]^2\right)
\end{equation}
a.s., where
\begin{align*}
I_1(\tau)&= \int_0^\tau \frac{d\sigma}{(\tau-\sigma)^{\frac{1}{2}+\alpha}}
\left(\int_\sigma^\tau d\rho \Vert z_M(.,\rho)\Vert_2^2\right)^{\frac{1}{2}},\\
I_2(\tau)&=\int_0^\tau \frac{d\sigma}{(\tau-\sigma)^{-\frac{\delta}{2}+\alpha+1}}
\left(\int_0^\sigma d\rho (\sigma-\rho)^{-\delta}\Vert z_M(.,\rho)\Vert_2^2\right)^{\frac{1}{2}},\\
I_3(\tau)&= \int_0^\tau \frac{d\sigma }{(\tau-\sigma)^{\alpha+1}}
\int_\sigma^\tau d\rho\left(\frac{1}{(\rho-\sigma)^\alpha}+\frac{1}{(\tau-\rho)^\alpha}\right )\Vert z_M(.,\rho)\Vert_2,\\
I_4(\tau)&=\int_0^\tau \frac{d\sigma }{(\tau-\sigma)^{\alpha+1}}
\left(\int_\sigma^\tau d\rho \int_\sigma^\rho d\xi \frac{\Vert z_M(.,\rho)-z_M(.,\xi)\Vert_2}{(\rho-\xi)^{\alpha+1}}\right),\\
I_5(\tau)&=\int_0^\tau \frac{d\sigma}{(\tau-\sigma)^{-\frac{\delta}{2}+\alpha+1}}
\int_0^\sigma \frac{d\rho}{(\sigma-\rho)^{\frac{\delta}{2}}}\left(\frac{1}{\rho^\alpha}+\frac{1}{(\sigma-\rho)^\alpha}\right) \Vert z_M(.,\rho)\Vert_2,\\
I_6(\tau)&=\int_0^\tau \frac{d\sigma}{(\tau-\sigma)^{-\frac{\delta}{2}+\alpha+1}}
\int_0^\sigma \frac{d\rho}{(\sigma-\rho)^{\frac{\delta}{2}}} \int_0^\rho d\xi \frac{\Vert z_M(.,\rho)-z_M(.,\xi)\Vert_2}{(\rho-\xi)^{\alpha+1}}.\\
\end{align*}
Set $T_k(t)= \int_0^t d\tau\ \left[I_k(\tau)\right]^2$, $k=1, \ldots, 6$.

The function $\sigma\mapsto (\tau-\sigma)^{-\frac{1}{2}-\alpha}$ is integrable on $(0,\tau)$ for $\alpha\in(0,\frac{1}{2})$. Thus
we have 
\begin{equation} 
\label{316}
T_1(t)\le c \int_0^t d\tau  \Vert z_M(.,\tau)\Vert_2^2.
\end{equation}
 Since we can choose $\delta>2\alpha$, we have that $\sigma\mapsto (\tau-\sigma)^{-\alpha-1+\frac{\delta}{2}}$ is integrable on $(0,\tau)$. Then, applying Schwarz inequality with respect to the measure 
 given by $d\sigma(\tau-\sigma)^{-\alpha-1+\frac{\delta}{2}}$,  we obtain
 \begin{align}
 \label{317}
T_2(t)&\le c \int_0^t d\tau\ \int_0^\tau d\sigma\ (\tau-\sigma)^{-\alpha-1+\frac{\delta}{2}}\left(\int_0^\sigma d\rho\ (\sigma-\rho)^{-\delta}
\Vert  z_M(.,\rho)\Vert_2^2\right)\nonumber\\
&\le c\int_0^t d\tau \left(\sup_{0\le \rho\le \tau}\Vert z_M(\cdot,\rho\Vert_2)^2\right)\left(\int_0^\tau d\sigma (\tau-\sigma)^{-\alpha-1+\frac{\delta}{2}}\sigma^{1-\delta}\right)\nonumber\\
&\le c\int_0^t d\tau \left( \sup_{0\le \rho\le \tau}\Vert z_M(\cdot,\rho\Vert_2)^2\right),
\end{align}
where in the last inequality we have used that $\alpha+\frac{\delta}{2}<1$ along with the definition of Euler's Beta function.

By integrating one obtains
\begin{equation*}
\int_\sigma^\tau d\rho\left(\frac{1}{(\rho-\sigma)^\alpha}+\frac{1}{(\tau-\rho)^\alpha}\right )= \frac{2(\tau-\sigma)^{1-\alpha}}{1-\alpha}.
\end{equation*}
Moreover, the function $\sigma\mapsto (\tau-\sigma)^{-2\alpha}$ is integrable on $(0,\tau)$. Consequently,
\begin{align}
\label{318}
T_3(t)&\le c \int_0^t d\tau \left(\sup_{\rho\in[0,\tau]} \Vert  z_M(.,\rho)\Vert_2^2\right) \int_0^\tau d\sigma\ (\tau-\sigma)^{-2\alpha}\nonumber\\
 &\le c  \int_0^t d\tau \left( \sup_{\rho\in[0,\tau]} \Vert  z_M(.,\rho)\Vert_2^2\right).
\end{align}
For any $\tau\in (0,t)$, set
\begin{equation*}
I_\tau= \int_0^\tau d\sigma \ (\tau-\sigma)^{-\alpha-1+\frac{\delta}{2}}\left(\int_0^\sigma d\rho\ (\sigma-\rho)^{-\frac{\delta}{2}}\left(\frac{1}{\rho^\alpha}
+\frac{1}{(\sigma-\rho)^\alpha}\right)\right).
\end{equation*}
It is a simple exercise to check that for $\alpha+\frac{\delta}{2}<1$,
$\sup_{\tau\in[0,t]} I_\tau < +\infty$.

Since
\begin{equation*}
T_5(t) \le \int_0^t d\tau\ I_\tau^2 \left(\sup_{\rho\in[0,\tau]} \Vert  z_M(.,\rho)\Vert_2^2\right),
\end{equation*}
we conclude that
\begin{equation}
\label{319}
T_5(t) \le c \int_0^t d\tau \left(\sup_{\rho\in[0,\tau]} \Vert  z_M(.,\rho)\Vert_2^2\right).
\end{equation}

Fix $\eta\in(0,1)$ so that $\sigma\mapsto (\tau-\sigma)^{-\eta}$ is integrable on $(0,\tau)$. Applying Schwarz inequality first with respect to the
measure $d\sigma (\tau-\sigma)^{-\eta}$, and then with respect to the Lebesgue measure on the interval $(\sigma,\tau)$ yields
\begin{align*}
T_4(t)&= \int_0^t d\tau \Big( \int_0^\tau \frac{d\sigma}{(\tau-\sigma)^\eta} (\tau-\sigma)^{-\alpha-1+\eta}\\
&\quad \times \Big(\int_\sigma^\tau d\rho \int_\sigma^\rho d\xi \frac{\Vert z_M(.,\rho)-z_M(.,\xi)\Vert_2}{(\rho-\xi)^{\alpha+1}}\Big)\Big)^2\\
&\le c \int_0^t d\tau\  \int_0^\tau \frac{d\sigma}{(\tau-\sigma)^\eta} (\tau-\sigma)^{-2\alpha-2+2\eta}\\
&\quad \times \left(\int_\sigma^\tau d\rho \int_\sigma^\rho d\xi \frac{\Vert z_M(.,\rho)-z_M(.,\xi)\Vert_2}{(\rho-\xi)^{\alpha+1}}\right)^2\\
&\le c \int_0^t d\tau\  \int_0^\tau d\sigma\ (\tau-\sigma)^{\eta -2\alpha-1}\\
&\quad \times  \int_\sigma^\tau d\rho \left( \int_\sigma^\rho d\xi \frac{\Vert z_M(.,\rho)-z_M(.,\xi)\Vert_2}{(\rho-\xi)^{\alpha+1}}\right)^2.
\end{align*}
By choosing $\eta>2\alpha$, the function $\sigma\mapsto (\tau-\sigma)^{\eta -2\alpha-1}$ is integrable on $(0,\tau)$. Thus, from the preceding inequalities
we obtain
\begin{align}
\label{320}
T_4(t)&\le c \int_0^t d\tau\  \int_0^\tau d\rho \left( \int_0^\rho d\xi \frac{\Vert z_M(.,\rho)-z_M(.,\xi)\Vert_2}{(\rho-\xi)^{\alpha+1}}\right)^2\nonumber\\
&\le c \int_0^t d\tau\  \Vert z_M\Vert_{\alpha,2,\tau}^2\ .
\end{align}
By Fubini's theorem and evaluations based upon Euler's Beta function, we have
\begin{align}
\label{321}
T_6(t)&= \int_0^t d\tau \Big(\int_0^\tau d\rho \Big( \int_\rho^\tau d\sigma (\tau-\sigma)^{\frac{\delta}{2}-\alpha-1}(\sigma-\rho)^{-\frac{\delta}{2}}\Big)\nonumber\\
&\quad\times \int_0^\rho d\xi \frac{\Vert z_M(.,\rho)-z_M(.,\xi)\Vert_2}{(\rho-\xi)^{\alpha+1}}\Big)^2\nonumber\\
&\le c \int_0^t d\tau\left(\int_0^\tau d\rho \left(\int_0^\rho d\xi \frac{\Vert z_M(.,\rho)-z_M(.,\xi)\Vert_2}{(\rho-\xi)^{\alpha+1}}\right)^2\right)\nonumber\\
&\le c \int_0^t d\tau\  \Vert z_M\Vert_{\alpha,2,\tau}^2.
\end{align}


Finally, inequalities (\ref{3160}) to (\ref{321}) imply
\begin{equation}
\label{322}
\Vert z_M\Vert_{\alpha,2,t}^2 \le R \int_0^t d\tau\ \Vert z_M\Vert_{\alpha,2,\tau}^2\ 
\end{equation}
a.s.
By Gronwall's lemma, this clearly implies the uniqueness of the mild solution.
Now the proof of  part (b) of Theorem \ref{t1} is complete. \hfill $\blacksquare$
\bigskip

\noindent{\bf Proof of Statement $(c)$ of Theorem \ref{t1}}
\bigskip

We investigate
each of the functions (\ref{19})--(\ref{21}) separately.

\begin{proposition}
\label{p2}
Assume that Hypotheses ($K_{\beta,\beta^{\prime}}$) and ($I$) hold. Then, there exists $c\in(0,+\infty)$ 
such that the estimate
\begin{equation}
\left\|A(\varphi)(.,t) -A(\varphi)(.,s)\right\|  _{2}\leq
c\left|  t-s\right|  ^{\theta^{\prime}} \label{33}
\end{equation}
\textit{holds for all }$s,t\in\left[  0,T\right]  $ and every $\theta^{\prime
}\in\left(  0,\frac{\beta}{2}\right]  $.
\end{proposition}

\noindent\textbf{Proof. }Relation (\ref{19}) defines a classical solution to
(\ref{2}) when $g=h=0$, so that the standard regularity theory for
linear parabolic equations gives $(x,t)\mapsto A(\varphi)
(x,t)\in\mathcal{C}^{\beta,\frac{\beta}{2}}(\overline{D}\times\left[
0,T\right]  \mathbb{)}$ (see, for instance, \cite{eidelzhita}), from which
(\ref{33}) follows immediately. \hfill $\blacksquare$

\bigskip

Regarding (\ref{20}) we have the following result.
\begin{proposition}
\label{p3}
Assume that the same hypotheses as in 
Theorem \ref{t1} (a) hold and let $u_{M}$ be any mild solution to (\ref{2}).
Then, there exists $c\in(0,+\infty)$ such that the
estimate
\begin{equation}
\left\|B(u_{M})(.,t) - B(u_{M})(.,s)\right\|  _{2}\leq
c\left|  t-s\right|  ^{\theta^{\prime\prime}}\left(1+\sup_{t\in[0,T]}\left\|  u_{M}(.,t)\right\|_2\right)  \label{34}
\end{equation}
holds a.s. for all $s,t\in\left[  0,T\right]$ and every
$\theta^{\prime\prime}\in\left(  0,\frac{1}{2}\right)$.
\end{proposition}

\noindent\textbf{Proof. } Without restricting the generality, we may assume that $t>s$.
We have
\begin{align}
&B(u_{M})(.,t) - B(u_{M})(.,s)
=\int_{s}^{t}d\tau\int_{D}dy\ G(.,t;y,\tau)g\left(  u_{M}(y,\tau)\right)
\nonumber\\
&\quad \quad  +\int_{0}^{s}d\tau\int_{D}dy\left(  G(.,t;y,\tau)-G(.,s;y,\tau)\right)
g\left(  u_{M}(y,\tau)\right),  \label{35}
\end{align}
and remark that in order to keep track of the increment $t-s$ we can estimate
the first term on the right-hand side of (\ref{35}) by using the same
kind of arguments as we did in the first part of the proof of Lemma \ref{l2}. For
every $x\in D$ this gives
\begin{align*}
&  \int_{s}^{t}d\tau\int_{D}dy\ \left|  G(x,t;y,\tau)g\left(  u_{M}
(y,\tau)\right)  \right| \\
&  \leq c(t-s)^{\frac{1}{2}}\left(  \int_{s}^{t}d\tau\int_{D}dy\left|
G(x,t;y,\tau)\right|  \left(  1+\left|  u_{M}(y,\tau)\right|  ^{2}\right)
\right)  ^{\frac{1}{2}},
\end{align*}
so that we eventually obtain
\begin{equation}
\left\|  \int_{s}^{t}d\tau\int_{D}dy\ G(.,t;y,\tau)g\left(  u_{M}(y,\tau)\right)  \right\|  _{2}
\leq c(t-s)^{\frac{1}{2}}\left(  1+\sup_{t\in[0,T]}\left\| u_{M}(.,t)\right\|_2\right)
\label{36}
\end{equation}
a.s. for all $s,t\in\left[  0,T\right]  $ with $t>s$. In \ a similar manner,
we can keep track of the increment $t-s$ in the second term on the right-hand
side of (\ref{35}) by using (\ref{15}). We thus have
\begin{align}
&  \left\|  \int_{0}^{s}d\tau\int_{D}dy\left(  G(.,t;y,\tau)-G(.,s;y,\tau
)\right)  g\left(  u_{M}(y,\tau)\right)  \right\|  _{2}^{2}\nonumber\\
&  \leq c\int_{0}^{s}d\tau\int_{D}dy\int_{D}dx\left|  G(x,t;y,\tau
)-G(x,s;y,\tau)\right|  \left(  1+\left|  u_{M}(y,\tau)\right|  ^{2}\right)\nonumber \\
&  \leq c(t-s)^{\delta}\int_{0}^{s}d\tau(s-\tau)^{-\delta}\left(  1+\left\|
u_{M}(.,\tau)\right\|  _{2}^{2}\right)\nonumber \\
&  \leq c(t-s)^{\delta}\left(  1+\sup_{t\in[0,T]}\left\|u_{M}(.,t)\right\|_2^{2}\right)\label{37}
\end{align}
for every $\delta\in\left(  \frac{d}{d+2},1\right)  $, 
a.s. for all $s,t\in\left[  0,T\right]  $ with $t>s$.  This last relation holds
\textit{a fortiori} for each $\delta\in\left(  0,1\right)$, so that
(\ref{36}) and (\ref{37}) indeed prove (\ref{34}). 

 \hfill $\blacksquare$

\bigskip

As for the stochastic term (\ref{21}), we have the following.

\begin{proposition}
\label{p4}
Assume the same hypotheses as in 
Theorem \ref{t1} (c), and let $u_{M}$ be any mild solution to (\ref{2}).
Then, there exists $c\in(0,+\infty)$ such that the
estimate
\begin{equation}
\left\|C(u_{M})(.,t) - C(u_{M})(.,s)\right\|  _{2}\leq
cr_{\alpha}^{H}\left|  t-s\right|  ^{\theta^{\prime
\prime\prime}}\left(  1+\left\|  u_{M}\right\|_{\alpha,2,T}\right)
\label{38}
\end{equation}
holds a.s. for all $s,t\in\left[  0,T\right]$ and every
$\theta^{\prime\prime\prime}\in\left(  0,\frac{1}{2}-\alpha\right)$.
\end{proposition}

The proof of Proposition \ref{p4} is more complicated than that of Proposition \ref{p3}.
We begin with a
preparatory result whose proof is based on inequalities
(\ref{15})--(\ref{17}).


\begin{lemma}
\label{l3}
With the same hypotheses as in part (a) of Theorem \ref{t1}, 
 the estimates
\begin{equation}
\sup_{i\in \mathbb{N}^+}\left\|  f_{i,t,s}^{\ast}(u_{M})(.,\tau)\right\|  _{2}\leq c\left(
t-s\right)  ^{\frac{\delta}{2}}\left(  s-\tau\right)  ^{-\frac{\delta}{2}
}\left(1+\sup_{t\in[0,T]}\left\| u_{M}(.,t)\right\|_2\right)  \label{40}
\end{equation}
\textit{and }
\begin{align}
& \sup_{i\in\mathbb{N}^+}\left\|  f_{i,t,s}^{\ast}(u_{M})(.,\tau)-f_{i,t,s}^{\ast}(u_{M}
)(.,\sigma)\right\|  _{2}\nonumber\\
&  \leq c\left(  t-s\right)  ^{\frac{\delta}{2}}\left(  s-\tau\right)
^{-\frac{\delta}{2}}\left\|  u_{M}(.,\tau)-u_{M}(.,\sigma)\right\|
_{2}\nonumber\\
&  +c\left(  t-s\right)  ^{\frac{\delta}{2}}\left(  s-\tau\right)  ^{-\frac
{1}{2}}\left(  \tau-\sigma\right)  ^{\frac{1}{2}(1-\delta)}\left(1+\sup_{t\in[0,T]}\left\|u_{M}(.,t)\right\|_2\right)  \label{41}
\end{align}
hold a.s. for every $\delta\in\left(  \frac{d}{d+2},1\right)$ and for all $\sigma,\tau\in\left[  0,s\right)$ with $\tau>\sigma$ in
(\ref{41}).
\end{lemma}

\noindent\textbf{Proof. }The proof of (\ref{40}) is analogous to that of
(\ref{25}) and is thereby omitted. As for (\ref{41}), by using
 Schwarz inequality relative to the measures $dy\left|
G(x,t;y,\tau)-G(x,s;y,\tau)\right|$ and
\[
dy\left|  G(x,t;y,\tau)-G(x,s;y,\tau)-G(x,t;y,\sigma)+G(x,s;y,\sigma)\right|
\]
on $D$ along with Hypothesis (L) for $h$, we get
\begin{align}
&  \left\|  f_{i,t,s}^{\ast}(u_{M})(.,\tau)-f_{i,t,s}^{\ast}(u_{M}
)(.,\sigma)\right\|  _{2}^{2}\nonumber\\
&  \leq c\int_{D}dx\int_{D}dy\left|  G(x,t;y,\tau)-G(x,s;y,\tau)\right|
\left|  u_{M}(y,\tau)-u_{M}(y,\sigma)\right|  ^{2}\nonumber\\
&  +c\int_{D}dx\int_{D}dy\left|  G(x,t;y,\tau)-G(x,s;y,\tau)-G(x,t;y,\sigma
)+G(x,s;y,\sigma)\right|\nonumber\\
&\quad\quad \times \left(  1+\left|  u_{M}(y,\sigma)\right|
^{2}\right) \nonumber\\
&  \leq c\left(  t-s\right)  ^{\delta}\left(  s-\tau\right)  ^{-\delta
}\left\|  u_{M}(.,\tau)-u_{M}(.,\sigma)\right\|  _{2}^{2}\nonumber\\
&\leq c\left(  t-s\right)  ^{\delta}\left(  s-\tau\right)  ^{-1}\left(
\tau-\sigma\right)  ^{1-\delta}\left(1+\sup_{t\in[0,T]}\left\|u_{M}(.,t)\right\|_2^2\right),  \label{42}
\end{align}
a.s. for all $s,t,\sigma,\tau\in\left[  0,T\right]  $ with $t\geq
s>\tau>\sigma$ and every $\delta\in\left(  \frac{d}{d+2},1\right)  $, as a
consequence of (\ref{15}), (\ref{44}) and the Gaussian property.

\hfill $\blacksquare$
\bigskip

\noindent\textbf{Proof of Proposition \ref{p4}.} For $t>s$ we write
\begin{align}
C(u_{M})(.,t) - C(u_{M})(.,s) &= \sum_{i=1}^{+\infty}\lambda_{i}^{\frac{1}{2}}\int_{s}^{t}f_{i,t}
(u_{M})(.,\tau)B_{i}^{H}(d\tau)\nonumber\\
& +\sum_{i=1}^{+\infty}\lambda_{i}^{\frac{1}{2}}\int_{0}^{s}f_{i,t,s}^{\ast
}(u_{M})(.,\tau)B_{i}^{H}(d\tau). \label{45}
\end{align}
In order to estimate the first
term on the right-hand side of
(\ref{45}), we can start by using inequalities (\ref{24}) and
(\ref{25}) to obtain
\begin{align}
&  \sum_{i=1}^{+\infty}\lambda_{i}^{\frac{1}{2}}\left\|  \int_{s}^{t}
f_{i,t}(u_{M})(.,\tau)B_{i}^{H}(d\tau)\right\|  _{2}
\nonumber\\
&  \leq cr_{\alpha}^{H}\left(  \int_{s}^{t}\frac{d\tau}
{(\tau-s)^{\alpha}}+\int_{s}^{t}d\tau\frac{\left\|  u_{M}(.,\tau)\right\|
_{2}^{{}}}{(\tau-s)^{\alpha}}+\int_{s}^{t}d\tau\int_{s}^{\tau}d\sigma
\frac{\left\|  u_{M}(.,\tau)-u_{M}(.,\sigma)\right\|  _{2}^{{}}}{\left(
\tau-\sigma\right)  ^{\alpha+1}}\right. \nonumber\\
&  +\left.  \int_{s}^{t}d\tau(t-\tau)^{-\frac{\delta}{2}}\int_{s}^{\tau
}d\sigma(\tau-\sigma)^{\frac{\delta}{2}-\alpha-1}\left(  1+\left\|
u_{M}(.,\sigma)\right\|  _{2}\right)  \right)  \label{46}
\end{align}
a.s. for every $s,t$ $\in\left[  0,T\right]  $ with $t>s$ and each $\delta
\in\left(  \frac{d}{d+2},1\right)$.

Furthermore, we have
\begin{align}
&  \int_{s}^{t}\frac{d\tau}{(\tau-s)^{\alpha}}+\int_{s}^{t}d\tau\frac{\left\|
u_{M}(.,\tau)\right\|  _{2}^{{}}}{(\tau-s)^{\alpha}}+\int_{s}^{t}d\tau\int
_{s}^{\tau}d\sigma\frac{\left\|  u_{M}(.,\tau)-u_{M}(.,\sigma)\right\|
_{2}^{{}}}{\left(  \tau-\sigma\right)  ^{\alpha+1}}\nonumber\\
&  \leq c\left(  \left(  t-s\right)  ^{1-\alpha}\left(  1+\left\|
u_{M}\right\|  _{\alpha,2,T}\right)  +(t-s)^{\frac{1}{2}}\left\|
u_{M}\right\|  _{\alpha,2,T}\right) \nonumber\\
&  \leq c(t-s)^{\frac{1}{2}}\left(  1+\left\|  u_{M}\right\|_{\alpha,2,T}\right)  \label{47}
\end{align}
since $\alpha<\frac{1}{2}$. Moreover,
\begin{align}
&  \int_{s}^{t}d\tau(t-\tau)^{-\frac{\delta}{2}}\int_{s}^{\tau}d\sigma
(\tau-\sigma)^{\frac{\delta}{2}-\alpha-1}\left(  1+\left\|  u_{M}
(.,\sigma)\right\|  _{2}\right) \nonumber\\
&  \le c\ \left(1+\sup_{t\in[0,T]}\left\|  u_{M}(.,t)\right\|_2\right)\int_{s}^{t}d\tau\  (t-\tau)^{-\frac{\delta}{2}}\int_s^\tau d\sigma\ \left(\tau-\sigma\right)
^{\frac{\delta}{2}-\alpha-1} \nonumber\\
&  \leq c(t-s)^{1-\alpha}\left(1+\sup_{t\in[0,T]}\left\|  u_{M}(.,t)\right\|_2\right),  \label{48}
\end{align}
by virtue of the convergence of the
integral, which can be expressed in terms of Euler's Beta function 
since $\alpha<\frac{\delta}{2}$. The substitution of (\ref{47}) and (\ref{48}) into
(\ref{46}) then leads to the inequality
\begin{equation}
\sum_{i=1}^{+\infty}\lambda_{i}^{\frac{1}{2}}\left\|  \int_{s}^{t}
f_{i,t}(u_{M})(.,\tau)B_{i}^{H}(d\tau)\right\|  _{2}\leq
cr_{\alpha}^{H}(t-s)^{\frac{1}{2}}\left(  1+\left\|
u_{M}\right\|_{\alpha,2,T}\right)  \label{hoelder8}
\end{equation}
a.s. for every $s,t$ $\in\left[  0,T\right]  $ with $t>s$. 

It remains to
estimate the second term on the right-hand side of (\ref{45}). 
From (\ref{i}) with $T$ replaced by $s$, we have
\begin{align*}
&  \sum_{i=1}^{+\infty}\lambda_{i}^{\frac{1}{2}}\left\|  \int_{0}^{s}
f_{i,t,s}^{\ast}(u_{M})(.,\tau)B_{i}^{H}(d\tau)\right\|
_{2}
\leq r_\alpha^H\\ 
&\times\sup_{i\in\mathbb{N}^+}\int_{0}^{s}d\tau\left(  \frac{\left\|
f_{i,t,s}^{\ast}(u_{M})(.,\tau)\right\|  _{2}}{\tau^{\alpha}}+\int
_{0}^{\tau}d\sigma\frac{\left\|  f_{i,t,s}^{\ast}(u_{M})(.,\tau)-f_{i,t,s}
^{\ast}(u_{M})(.,\sigma)\right\|  _{2}}{\left(  \tau-\sigma\right)
^{\alpha+1}}\right).
\end{align*}

By substituting (\ref{40}) and (\ref{41}) we  obtain
\begin{align}
&  \sum_{i=1}^{+\infty}\lambda_{i}^{\frac{1}{2}}\left\|  \int_{0}^{s}
f_{i,t,s}^{\ast}(u_{M})(.,\tau)B_{i}^{H}(d\tau)\right\|
_{2}\nonumber\\
&  \leq cr_{\alpha}^{H}(t-s)^{\frac{\delta}{2}}\left(
\int_{0}^{s}d\tau\left(  s-\tau\right)  ^{-\frac{\delta}{2}}\tau^{-\alpha
}\left(1+\sup_{t\in[0,T]}\left\|  u_{M}(.,t)\right\|_2\right)  \right.
\nonumber\\
&  +\int_{0}^{s}d\tau\left(  s-\tau\right)  ^{-\frac{\delta}{2}}\int_{0}
^{\tau}d\sigma\frac{\left\|  u_{M}(.,\tau)-u_{M}(.,\sigma)\right\|_{2}
}{\left(  \tau-\sigma\right)  ^{\alpha+1}}\nonumber\\
&  +\left.  \int_{0}^{s}d\tau\left(  s-\tau\right)  ^{-\frac{1}{2}}\int
_{0}^{\tau}d\sigma\left(  \tau-\sigma\right)  ^{\frac{1}{2}\left(
1-\delta\right)  -\alpha-1}\left(1+\sup_{t\in[0,T]}\left\|  u_{M}(.,t)\right\|_2\right)  \right) \nonumber\\
&  \leq c\ r_{\alpha}^{H}(t-s)^{\frac{\delta}{2}}\left(1+\left\|  u_{M}\right\|_{\alpha,2,T}\right)\nonumber\\
&\quad\times \left(1+\int_{0}^{s}d\tau\left(  s-\tau\right)  ^{-\frac{1}{2}}\int_{0}^{\tau}d\sigma
\left(  \tau-\sigma\right)  ^{\frac{1}{2}\left(  1-\delta\right)-\alpha-1}\right),\label{49}
\end{align}
where we have got the last estimate using Schwarz inequality with
respect to the measure $d\tau$ on $\left(  0,s\right)$ along with
(\ref{3}) in the first two integrals on the right-hand side.

By imposing the additional restriction $\delta<1-2\alpha$, we have
\begin{equation*}
\int_{0}^{s}d\tau\left(  s-\tau\right)  ^{-\frac{1}{2}}\int_{0}^{\tau}d\sigma
\left(  \tau-\sigma\right)  ^{\frac{1}{2}\left(  1-\delta\right)-\alpha-1}
<+\infty.
\end{equation*}
Thus, we have proved that
\begin{equation}
\sum_{i=1}^{+\infty}\lambda_{i}^{\frac{1}{2}}\left\|  \int_{0}^{s}
f_{i,t,s}^{\ast}(u_{M})(.,\tau)B_{i}^{H}(d\tau)\right\|
_{2}\leq cr_{\alpha}^{H}(t-s)^{\frac{\delta}{2}}\left(
1+\left\|  u_{M}\right\|_{\alpha,2,T}\right)  \label{50}
\end{equation}
a.s. for all $s,t$ $\in\left[  0,T\right]  $ with $t>s$ and every $\delta
\in\left( \frac{d}{d+2},1-2\alpha\right)$. The existence of this restricted
interval of values of $\delta$ is possible by our choice of $\alpha$. Relations
(\ref{45}), (\ref{hoelder8}) and (\ref{50}) clearly 
yield (\ref{38}) with 
$\theta^{\prime\prime\prime}
=\frac{\delta}{2}\in\left(  0,\frac{1}{2}-\alpha\right)  $.
 \hfill $\blacksquare$

\bigskip

It is immediate that Propositions \ref{p2} to \ref{p4} imply statement $(c)$ of  Theorem \ref{t1}.
Notice that $R_\alpha^H = c(1+r_\alpha^H)$, with $r_\alpha^H$ defined in (\ref{29}).
\bigskip

Finally, we give an alternate to the result proved before, as mentioned in Section \ref{s2}, Remark 5.

\begin{proposition}
\label{pfact}
The assumptions are as in Theorem \ref{t1} part (a). Then,
\begin{equation}
\Vert C(u_{M})(.,t)-C(u_{M})(.,s)\Vert_{2}\leq R |t-s|^{\theta^{\prime\prime\prime\prime}}
\left(1+\Vert u_{M}\Vert_{\alpha,2,T}\right)
\label{323}
\end{equation}
holds a.s. for all $s,t\in\lbrack0,T]$ and every 
$\theta^{\prime\prime\prime\prime}\in\left(  0,\frac{2}{d+2}\wedge\frac{1}{2}\right)$. Consequently, 
\beq
\label{3231}
\Vert u_{M})(.,t)-u_{M})(.,s)\Vert_2\le R |t-s|^\theta \left(1+\Vert u_{M}\Vert_{\alpha,2,T}\right),
\eeq
a.s. for all $s,t\in\lbrack0,T]$ and each $\theta\in\left(0,\frac{2}{d+2}\wedge\frac{\beta}{2}\right)$.
\end{proposition}
{\bf Proof.} 
We use the factorization method we alluded to in Section \ref{s2}. For this we express 
$C(u_{M})(.,t)$ in terms of the auxiliary $L^{2}(D)$-valued process
\begin{equation*}
Y_{\varepsilon}(u_{M})(.,t):=\sum_{i=1}^{+\infty}\lambda_{i}^{\frac{1}{2}}%
\int_{0}^{t}(t-\tau)^{-\varepsilon}f_{i,t}(u_{M})(.,\tau)B_{i}^H(d\tau)
\end{equation*}
defined for every $\varepsilon\in\left(  0,\frac{1}{2}\right)  $. In fact, by
repeated applications of Fubini's theorem and by using the fundamental
property $U(t,\tau)U(\tau,\sigma)=U(t,\sigma)$ for the evolution operators
defined in (\ref{30}) we obtain
\begin{align}
C(u_{M})(.,t)  &  =\sum_{i=1}^{+\infty}\lambda_{i}^{\frac{1}{2}}\int_{0}%
^{t}f_{i,t}(u_{M})(.,\tau)B_{i}^{H}(d\tau)\nonumber\\
&  =\frac{\sin(\varepsilon\pi)}{\pi}\int_{0}^{t}d\tau(t-\tau)^{\varepsilon
-1}\int_{D}dyG(.,t;y,\tau)Y_{\varepsilon}(u_{M})(y,\tau) \label{factorization}%
\end{align}
for every $t\in\left[  0,T\right]$ a.s. 

Next we prove that a.s., 
\begin{equation}
\sup_{t\in\lbrack0,T]}\Vert Y_{\varepsilon}(u_{M})(.,t)\Vert_{2}\leq R\left(
1+\Vert u_{M}\Vert_{\alpha,2,T}\right).  \label{324}
\end{equation}
Indeed, according with (\ref{i}),
\begin{align*}
&\Vert Y_{\varepsilon}(u_{M})(.,t)\Vert_{2}\le r_\alpha^H \sup_{i\in\mathbb{N}^+} \int_0^t d\tau (t-\tau)^{-\varepsilon}\\
&\quad\times\left\{\frac{\Vert f_{i,t}(u)(.,\tau)\Vert_2}{\tau^\alpha}+\int_0^\tau d\sigma \frac{\Vert f_{i,t}(u)(.,\tau)-f_{i,t}(u)(.,\sigma)\Vert_2 }{(\tau-\sigma)^{\alpha+1}}\right\}.
\end{align*}
From the estimate (\ref{24}) and the definition of the Beta function, we have
\begin{equation*}
\int_0^t d\tau (t-\tau)^{-\varepsilon}\frac{\Vert f_{i,t}(u)(.,\tau)\Vert_2}{\tau^\alpha}\le c\left(1+\sup_{\tau\in[0,t]} \Vert u(.,\tau)\Vert_2\right).
\end{equation*}
Since $\varepsilon\in(0,\frac{1}{2})$, applying Schwarz's inequality yields
\begin{align*}
&\int_0^t d\tau (t-\tau)^{-\varepsilon}\int_0^\tau \frac{\Vert u(.,\tau)-u(.,\sigma)\Vert_2}{(\tau-\sigma)^{\alpha+1}}\\
&\quad\le \left[\int_0^t  d\tau\ (t-\tau)^{-2\varepsilon}\right]^{\frac{1}{2}} \left[\int_0^t d\tau\ \left(\int_0^\tau d\sigma \frac{\Vert u(.,\tau)-u(.,\sigma)\Vert_2}{(\tau-\sigma)^{\alpha+1}}\right)^2\right]^{\frac{1}{2}}\\
&\quad\le c \Vert u\Vert_{\alpha,2,t}.
\end{align*}
Moreover, for any $\delta\in(0,\frac{1}{2})$ such that $\frac{\delta}{2}-\alpha>0$, we have
\begin{equation*}
\int_0^t d\tau (t-\tau)^{-\varepsilon-\frac{\delta}{2}} \int_0^\tau d\sigma (\tau-\sigma)^{\frac{\delta}{2}-\alpha-1}\left(1+\Vert u(.,\sigma)\Vert_2\right)
\le c\left(1+\sup_{\sigma\in[0,t]}\Vert u(.,\sigma)\Vert_2\right).
\end{equation*}
By virtue of (\ref{25}) the two above estimates imply
\begin{align*}
\int_0^t d\tau\int_0^\tau d\sigma (t-\tau)^{-\varepsilon}\frac{\Vert f_{i,t}(u)(.,\tau)-f_{i,t}(u)(.,\sigma)\Vert_2 }{(\tau-\sigma)^{\alpha+1}}
\le c \left(1+\Vert u\Vert_{\alpha,2,t}\right).
\end{align*}
This ends the  proof of (\ref{324}).
\smallskip

 We can now proceed by estimating the
time increments of $C(u_{M})$ using (\ref{factorization}) and (\ref{324}) rather than
with the expressions of Proposition \ref{p4}. For this
we  follow the arguments of the proof of (66) in Proposition 6 of \cite{sansovu1}
to see that,  by choosing
 $\theta^{\prime\prime\prime\prime}\in\left(0,\frac{2}{d+2}\wedge\frac{1}{2}\right)$ with the additional restriction
$\varepsilon\in\left(  \theta^{\prime\prime\prime\prime},\frac{2}{d+2}%
\wedge\frac{1}{2}\right)$, we obtain
\begin{align*}
&  \left\Vert C(u_{M})(.,t)-C(u_{M})(.,s)\right\Vert _{2}\\
&  \leq c\left(\left\Vert \int_{s}^{t}d\tau(t-\tau)^{\varepsilon-1}\int
_{D}dyG(.,t;y,\tau)Y_{\varepsilon}(u_{M})(y,\tau)\right\Vert _{2}\right. \\
&  +\left.  \left\Vert \int_{0}^{s}d\tau\int_{D}dy\left(  (t-\tau
)^{\varepsilon-1}G(.,t;y,\tau)-(s-\tau)^{\varepsilon-1}G(.,s;y,\tau)\right)
Y_{\varepsilon}(u_{M})(y,\tau)\right\Vert _{2}\right) \\
&  \leq R\left(  \left\vert t-s\right\vert ^{\varepsilon}+\left\vert
t-s\right\vert ^{\theta^{\prime\prime\prime\prime}}\right)  (1+\Vert
u_{M}\Vert_{\alpha,2,T})
\leq R |t-s|^{\theta^{\prime\prime\prime\prime}}\left(
1+\Vert u_{M}\Vert_{\alpha,2,T}\right),
\end{align*}
proving (\ref{323}).

Finally, this estimate along with those established in Propositions \ref{p2} and \ref{p3} provide (\ref{3231}) and finish the proof of
the proposition.

\hfill $\blacksquare$


\bigskip

\noindent\textbf{Acknowledgements.}
The research of the first author concerning this paper was completed at the Institute Mittag-Leffler
in Djursholm. 
The research of the second author was supported in part by the 
Institute of Mathematics of the University of Barcelona where this work was begun, and in
part by the ETH-Institute of Theoretical Physics in Zurich. They would like to
thank the three institutions for their very kind hospitality.

\end{document}